\documentclass[aop,MSNbibl,citesort,dvips]{arximspdf}

\doi{10.1214/10-AOP528}
\volume{38}
\issue{6}
\pubyear{2010}
\firstpage{2136}
\lastpage{2169}

\makeatletter

\let\comp\circ
\newcommand{\St}{\mathrm{St}}
\newcommand{\bigtimes}{\mathop{\mbox{\fontsize{16}{16}\selectfont$\times$}}_{j=1}^{m}}
\newtheorem{proposition}{Proposition}[section]
\newproclaim{definition}[proposition]{Definition}

\newproclaim{notation}[proposition]{Notation}
\newtheorem{theorem}[proposition]{Theorem}
\newtheorem{lemma}[proposition]{Lemma}

\newproclaim{remarks}[proposition]{Remarks}

\makeatother

\begin{document}
\begin{frontmatter}

\title{Multiple Stratonovich integral and Hu--Meyer formula for L{\'e}vy processes}
\runtitle{Multiple Stratonovich integral for L\'{e}vy processes}

\pdftitle{Multiple Stratonovich integral and Hu--Meyer formula for Levy processes}

\begin{aug}
\author[A]{\fnms{Merc{\`e}} \snm{Farr{\'e}}\ead[label=e1]{farre@mat.uab.cat}},
\author[A]{\fnms{Maria} \snm{Jolis}\thanksref{t1}\ead[label=e2]{mjolis@mat.uab.cat}} and
\author[A]{\fnms{Frederic} \snm{Utzet}\corref{}\thanksref{t1}\ead[label=e3]{utzet@mat.uab.cat}}
\runauthor{M. Farr{\'e}, M. Jolis and F. Utzet}
\affiliation{Universitat Aut{\`o}noma de Barcelona}
\address[A]{Department of Mathematics\\
Edifici C\\
Universitat Aut{\`o}noma de Barcelona\\
08193 Bellaterra (Barcelona)\\
Spain\\
\printead{e1}\\
\phantom{E-mail: }\printead*{e2}\\
\phantom{E-mail: }\printead*{e3}} 
\end{aug}

\thankstext{t1}{Supported by Grant MTM2009-08869 Ministerio de Ciencia
e Innovaci{\'o}n and FEDER.}

\pdfauthor{Merce Farre, Maria Jolis, Frederic Utzet}

\received{\smonth{3} \syear{2009}}
\revised{\smonth{11} \syear{2009}}

%
\begin{abstract}
In the framework of vector measures and the combinatorial approach to
stochastic multiple integral introduced by Rota and Wallstrom
[\textit{Ann. Probab.} \textbf{25} (1997) 1257--1283], we present an
It{\^o}
multiple integral and a Stratonovich multiple integral with respect
to a L\'{e}vy process with finite moments up to a convenient order. In
such a framework, the Stratonovich multiple integral is an integral
with respect to a product random measure whereas the It{\^o} multiple
integral corresponds to integrate with respect to a random measure
that gives zero mass to the diagonal sets. A general Hu--Meyer formula
that gives the relationship between both integrals is proved. As
particular cases, the classical Hu--Meyer formulas for the Brownian
motion and for the Poisson process are deduced. Furthermore, a
pathwise interpretation for the multiple integrals with respect to a
subordinator is given.
\end{abstract}

%

\pdfkeywords{60G51, 60H99, 05A18, 60G57, Levy processes,
Stratonovich integral, Hu--Meyer formula, random measures,
Teugels martingales}

\end{frontmatter}

\section{Introduction}\label{sec1}

Let $W=\{W_t, t\ge0\}$ be a standard Brownian motion. It{\^o}
\cite{ito-multiple} defined the multiple stochastic integral
of a function $f\in L^2({\mathbb{R}}_+^n,\break \mathcal{B}({\mathbb
{R}}_+^n), (dt)^{\otimes n})$,
\[
I_n(f)=\int\cdots\int_{{\mathbb{R}}^n_+} f(t_1,\ldots
,dt_n)\,dW_{t_1}\cdots dW_{t_n},
\]
taking care to ensure that the diagonal sets, like $\{(s_1,\ldots
,s_n)\in{\mathbb{R}}^n_+, s_1=s_2\}$,
do not contribute at all. For this reason the integral has very good
properties and is easy to work with. However,
for a function of the form
\[
(g_1\otimes\cdots\otimes g_n)(t_1,\ldots,t_n):=g(t_1)\cdots g(t_n),
\]
we have that, in general,
\[
I_n(g_1\otimes\cdots\otimes g_n)\ne I_1(g_1)\cdots I_1(g_n).
\]
That means the It\^{o} multiple integral does not behave like the
integral with respect to a product measure.

Many years later, Hu and Meyer \cite{hu-meyer1} introduced (although
they believed that this integral was already known
\cite{hu-meyer1}, page 75) a multiple integral, $I_n^S(f)$,
which followed the ordinary rules of multiple integration. They
called it the multiple Stratonovich integral.
Furthermore, Hu and Meyer
stated the relationship between the It{\^o} and Stratonovich
integrals, the celebrated Hu--Meyer formula,
adding the contribution of the diagonals to the It{\^o} integral: for
a function $f(t_1,\ldots,t_n)$ symmetric with good properties,
\[
I_n^S(f)
=\sum_{j=0}^{[n/2]}\frac{n!}{(n-2j)!j!2^j} I_{n-2j}
\biggl(\int_{{\mathbb{R}}_+^{j}}f(\bolds\cdot,t_1,t_1,t_2,t_2,\ldots
,t_{j},t_{j}) \,d t_{1}\cdots
d t_{j} \biggr).
\]
This formula is simple because the quadratic variation of the Brownian
motion is~$t$, and
the integral over coincidences of order three or superior are zero.
Following their ideas, Sol{\'e} and Utzet \cite{sole-pois1} proved a
Hu--Meyer formula for the Poisson process.
Again, in that case, the formula is relatively simple because the
variations of any order of the process can always
be written in terms of the Poisson process and $t$.

From another point of view, Engel \cite{engel}, working with a general
process with independent increments,
related the (It{\^o}) multiple stochastic integral with the theory of
vector valued measures, and Masani \cite{masani}, using also
vector valued measures and starting from the Wiener's original ideas,
developed both the It{\^o} and Stratonovich integrals (with respect to
the Brownian motion) and proved many profound results.
The vector measures approach is no simple matter; Engel's work covers
82 pages, and Masani's covers 160.
An important and clarifying contribution was made by Rota and Wallstrom
\cite{rota} who used combinatorial techniques
to show the features of the multiple stochastic integration. They did
not really work with integrals,
but with products
of vector measures. However, the path towards a general theory of
multiple stochastic integration had been laid.
See also P{\'e}rez--Abreu \cite{perez} for an interesting
generalization to Hilbert space valued random measures.
Further, Vershik and Tsilevich \cite{VerTsi03}, in a more algebraic context,
constructed a Fock factorization for a L\'{e}vy process, and some
important subspaces can be described through Rota and Wallstrom
concepts. We should also mention the very complete survey by Peccati
and Taqqu \cite{PecTaq08} in which a unified study
of multiple integrals, moments, cumulants and diagram formulas, as
well as applications to some new central limit theorems, is presented.

It is worth remarking that Rota and Walstrom's \cite{rota}
combinatorial approach to multiple integration
has been extended to the context of free probability in a very
interesting and fertile field of research, started by
Anshelevich (see \cite{Ans00,Ans02a,Ans02b,Ans04,Ans05} and the
references therein). In fact, Rota and Walstrom's ideas
fit very well with the combinatorics of free probability (see Nica and
Speicher \cite{NicSpe06}) and noncommutative L\'{e}vy
processes. Our renewed interest in Rota and Walstrom's paper
\cite{rota} was motivated by Anshelevich's work.

In the present paper we use the powerful Rota and Wallstrom's
\cite{rota} combinatorial machinery
to study the Stratonovich integral (the integral with respect to the
product random measure)
with respect to a L{\'e}vy processes with finite moments up to a
convenient order.
The key point is to understand how the product
of stochastic measures works on the diagonal sets, and that leads to
the diagonal measures defined by Rota and Wallstrom \cite{rota}.
For a L{\'e}vy process those measures are related
to the powers of the jumps of the process, and hence to a family of
martingales introduced by
Nualart and Schoutens \cite{nua-sch}, called Teugels martingales,
which offer excellent properties. Specifically,
these martingales
have deterministic predictable quadratic variation and this makes it
possible to easily construct
an It{\^o} multiple stochastic integral with respect to different
integrators, which can be interpreted as an integral with respect to a
random measure
that gives zero mass to the diagonal sets.
With all these ingredients we prove a general
Hu--Meyer formula. The paper uses arduous combinatorics because
of our need to work with stochastic multiple
integrals with respect to the different powers of the jumps of the
process, and such integrals
can be conveniently handled through the lattice of the partitions of a
finite set.

As in the Brownian case (see, e.g.,
\cite{johnson,sole-tra,hu-meyer2,masani}),
there are alternative methods to construct a multiple Stratonovich
integral based on approximation
procedures, and it is possible to relax the conditions on the
integrator process by
assuming more regularity on the integrand function. Such regularity is
usually expressed
in terms of the existence of \textit{traces} of the function
in a convenient sense.
The advantage of using
L{\'e}vy processes with finite moments lies in the fact that simple
$L^2(\Omega)$ estimates for the
multiple stochastic integral of simple functions can be obtained, and then
the multiple Stratonovich integral can be defined in an $L^2$ space
with respect to a measure 
that controls the behavior of the functions on the diagonal sets. In
this way, the problem of providing a manageable definition of the
traces is avoided.

We would like to comment that an impressive body of work on multiple
stochastic integrals with respect
to L{\'e}vy processes
has been done by Kallenberg, Kwapien, Krakowiak, Rosinski, Szulga,
Woyczinski and many others
(see \cite{kallenberg,krak,kwa,rosinski} and the references therein).
However, their approach is very different from ours, and assumes
different settings to those used in this work.
For this reason, we have only used
a few results by those authors.

The paper is organized as follows. In Section \ref{sec2} we review some
combinatorics concepts
and the basics of the stochastic measures as vector valued measures. In
Section \ref{random} we introduce the random measures
induced by a L{\'e}vy process, and we identify the diagonal measures in
such a case.
In Section \ref{measure} we study the relationship between the product
and It{\^o}
measures of a set, and we obtain
a Hu--Meyer formula for measures. In Section \ref{integral} we define
the multiple
It{\^o} stochastic integral and the multiple Stratonovich integral and also
prove the general Hu--Meyer formula for integrals. In Section \ref
{sec6}, as particular
cases, we deduce the classical Hu--Meyer formulas for the Brownian
motion and for the Poisson process.
We also study the case where the L\'{e}vy process is a subordinator,
and prove that both the
multiple It{\^o} stochastic integral and the multiple Stratonovich
integral can be computed in a pathwise sense.
Finally, in order to make the paper lighter, some of the combinatorial results
are included as an \hyperref[app]{Appendix}.

\section{Preliminaries}\label{sec2}

\subsection{Partitions of a finite set}
\label{partitions}
We need some notation of the combinatorics of the partitions of a
finite set; for details
we refer to Stanley \cite{stanley}, Chapter 3, or Rota and Wallstrom
\cite{rota}.

Let $F$ be a finite set.
A partition of $F$
is a family $\pi=\{B_1,\ldots, B_m\}$
of nonvoid subsets of $F$, pairwise disjoint,
such that $F=\bigcup_{i=1}^m B_i$. The elements $B_1,\ldots, B_m$ are called
the \textit{blocks} of the partition.
Denote by $\Pi(F)$ the set of all partitions of $F$, and write $\Pi_n$
for $\Pi(\{1,\ldots,n\})$.
Given $\sigma, \pi\in\Pi(F)$, we write $\sigma\le\pi$ if each
block of $\sigma$ is contained in some block of $\pi$; we then say
that $\sigma$ is a refinement of $\pi$. This relationship defines a
partial order that is called the \textit{reversed refinement order},
and it makes $\Pi(F)$ a lattice. We
write $\widehat0= \{\{x\}, x
\in F \}$, which is the minimal element, and
$\widehat1=\{F\}$ the maximal one.

We say that a partition $\pi\in\Pi(F)$ is of type
$(1^{r_1}2^{r_2}\cdots n^{r_n})$ if $\pi$ has exactly $r_1$ blocks
with 1
element, exactly $r_2$ blocks with 2 elements, and so on. In the same
way, for
$\sigma\le\pi, \# \sigma=m$ and $\# \pi=k$, we say that the segment
$[\sigma,\pi]$ is of type
$(1^{r_1}2^{r_2}\cdots m^{r_m})$ if there are exactly $r_1$ blocks of
$\pi$ in
$\sigma$; there are exactly $r_2$ blocks of $\pi$ that each one gives rise
to 2 blocks of $\sigma$, etc. Necessarily,
\[
\sum_{j=1}^m r_j=k \quad\mbox{and}\quad \sum_{j=1}^m j r_j=m.
\]
In that situation, the M{\"o}bius function of $[\sigma,\pi]$ is
\[
\mu(\sigma,\pi)=(-1)^{m-k}(2!)^{r_3}\cdots\bigl((m-1)!\bigr)^{r_m}.
\]
We use the M{\"o}bius inversion formula, that in the context of the lattice
of the partitions of a finite set, says that for two functions
$f,g\dvtx\Pi
(F)\longrightarrow\mathbb{R}$,
\[
g(\sigma)=\sum_{\pi\ge\sigma} f(\pi)\qquad \forall\sigma\in\Pi(F),
\]
if and only if
%
%
\begin{equation}
\label{moebius}
f(\sigma)=\sum_{\pi\ge\sigma}\mu(\sigma,\pi) g(\pi)\qquad \forall
\sigma\in
\Pi(F)
\end{equation}
(see \cite{stanley}, Proposition 3.7.2).

\subsection{Diagonal sets induced by a partition}
\label{subsets}
As we commented in the\break \hyperref[sec1]{Introduction}, we will introduce
two random measures on a $n$-dimensional space, and the diagonal sets
will play an essential role.
Diagonal sets can be conveniently described through the partitions of
the set
$\{1,\ldots,n\}$. We use the notation introduced by Rota and Wallstorm
\cite{rota}.

Let $S$ be an arbitrary set, and consider $C\subset S^n$.
Given $\pi\in\Pi_n$, we write $i \sim_\pi j$ if $i$ and $j$ belong
to the same block of $\pi$. Put
\[
{C} _{\ge\pi}=\{(s_1,\ldots, s_n)\in{C}\dvtx s_i=s_j \mbox{
if } i \sim_\pi j\}
\]
and
\[
{C} _{ \pi}=\{(s_1,\ldots, s_n)\in{C}\dvtx s_i=s_j \mbox{
if and only if } i \sim_\pi j\}.
\]
The sets $C_\pi$ are called \textit{diagonal sets}. Note that $C_\pi
=C\cap
S^n_\pi$
and $C_{\ge\pi}=C\cap S^n_{\ge\pi}$.

For example, for $n=4$ and
$\pi= \{\{1\},\{2\},\{3,4\} \}$, we have
\[
{C}_{\ge\pi}=\{(s_1,s_2,s_3,s_4)\in{C}\dvtx s_3=s_4\}
\]
and
\[
{C}_{\pi}=\{(s_1,s_2,s_3,s_4)\in{C}\dvtx s_3=s_4, s_1\ne s_2, s_1\ne
s_3, s_2\ne s_3\}.
\]

The sets corresponding to the minimal and maximal partitions are
specially important
\[
{C} _{\widehat0}=\{(s_1,\ldots, s_n)\in{C}\dvtx s_i\ne s_j, \forall
i\ne
j\}
\]
and
\[
{C}_{\widehat1}=\{(s_1,\ldots,
s_n)\in{C}\dvtx s_1= \cdots= s_n\}.
\]

If $\sigma\ne\pi$, then
%
%
\begin{equation}
\label{incompatibilitat}
C_\sigma\cap C_\pi=\varnothing
\quad\mbox{and}\quad (C_\pi)_\sigma=\varnothing.
\end{equation}
The above notation ${C} _{\ge\pi}$ is coherent with the reversed
refinement order
%
%
\begin{equation}
\label{union}
{C}_{\ge\pi} =\bigcup_{\sigma\ge\pi}{C} _ \sigma
\qquad(\mbox{disjoint union}).
\end{equation}
In particular,
$C=C_{\ge\widehat0}=\bigcup_{\sigma\in\Pi_n}{C} _ \sigma$.

\subsection{Random measures}\label{sec23}

Let $(\Omega, \mathcal{F}, {\mathbb{P}})$ be a complete probability space.
In this paper, a random measure $\Phi$ on a measurable space
$(S,\mathcal{S})$ is
an $L^2(\Omega)$-valued $\sigma$-additive vector measure, that means,
a map $\Phi\dvtx\mathcal{S}\to L^2(\Omega)$ such that
for every sequence $\{A_n, n\ge1\}\subset\mathcal{S}$, such that
$A_n\cap A_m=\varnothing, n\ne m$,
\[
\Phi\Biggl(\bigcup_{n=1}^\infty A_n\Biggr)=\sum_{n=1}^\infty\Phi(A_n)\qquad
\mbox{convergence in } L^2(\Omega).
\]

The $\sigma$-additive vector measures defined on a $\sigma$-field
inherit some basic properties of the ordinary measures, but not all.
So,
for a sake of easy reference, we write here a uniqueness property
translated to our setting. The proof is the same as the one for
ordinary measures.
\begin{proposition}
\label{monotona} Let $\Phi$ and $\Psi$ be two random measures on
$(S,\mathcal{S})$, and consider a family of sets $\mathcal{C}\subset
\mathcal{S}$
closed under finite intersection and such that $\sigma(\mathcal
{C})=\mathcal{S}$. Then
\[
\Phi=\Psi \mbox{ on } \mathcal{C} \quad\Longrightarrow\quad\Phi=\Psi
\mbox{ on }
\mathcal{S}.
\]
\end{proposition}

\subsection{\texorpdfstring{Product and It{\^o} stochastic measures}{Product and Ito stochastic measures}}\label{sec24}

Assume that the measurable space $(S,\mathcal{S})$ satisfies that for
every set $C\in\mathcal{S}^{\otimes n}$ and every
$\pi\in\Pi_n$, we
have $C_\pi\in\mathcal{S}^{\otimes n}$. As Rota and Wallstrom
\cite{rota} point out, this condition is satisfied if
$S$ is a Polish space and $\mathcal{S}$ its Borel $\sigma$-algebra.
We extend the definition of \textit{good random measure} introduced by
Rota and Wallstorm \cite{rota} to a family
of measures; specifically, we say that
the random measures $\Phi_1,\ldots, \Phi_k$ over a measurable space
$(S,\mathcal{S})$ are
\textit{jointly good random measures} if
the finite additive
product vector \textit{measure}
$\Phi_{1}\otimes\cdots\otimes\Phi_{k}$ defined on the product sets by
\[
(\Phi_{1}\otimes\cdots\otimes\Phi_{k} )(A_1\times\cdots\times
A_k)=\prod
_{j=1}^k\Phi_{i}(A_i),\qquad A_1,\ldots, A_k \in\mathcal{S},
\]
can be extended
to a (unique) $\sigma$-additive random measure on $(S^n,\mathcal
{S}^{\otimes n})$. This extension, obvious for
ordinary measures, is in general not transferred to arbitrary vector
measures (see Engel \cite{engel}, Masani \cite{masani} and Kwapien and
Woyczynski \cite{kwa}).

Given a good random measure $\Phi$ (in the sense that the $n$-fold product
$\Phi\otimes\cdots\otimes\Phi=\Phi^{\otimes n}$ satisfies the above
condition), the starting point of Rota
and Wallstrom~(\cite{rota}, Definition 1) is to consider new random
measures given by the restriction over the diagonal sets; specifically,
for $\pi\in\Pi_n$ they define
\[
\Phi_\pi^{\otimes n}(C):=\Phi^{\otimes n}(C_{\ge\pi})
\quad\mbox{and}\quad
\St^{[n]}_\pi(C):=\Phi^{\otimes n}(C_{\pi})\qquad
\mbox{for } C\in\mathcal{S}^{\otimes n}.
\]

The following definitions are the extension of these concepts
to a family of random
measures.
\begin{definition}Let $\Phi_{r_1},\ldots,\Phi_{r_n}$ be jointly good
random measures on $(S,\mathcal{S})$.
For a partition $\pi\in\Pi_n$,
define
%
%
\begin{equation}
\label{def-producte}
(\Phi_{r_1}\otimes\cdots\otimes\Phi_{r_n} )_\pi(C)=(\Phi
_{r_1}\otimes
\cdots\otimes\Phi_{r_n} )(C_{\ge\pi}),\qquad
C\in\mathcal{S}^{\otimes n},
\end{equation}
and
%
%
\begin{equation}
\label{def-ito}
\St^{(r_1,\ldots,r_n)}_\pi(C)=(\Phi_{r_1}\otimes\cdots\otimes\Phi_{r_n}
)(C_{\pi}),\qquad C\in\mathcal{S}^{\otimes n}.
\end{equation}
\end{definition}

In agreement with the notation in Rota
and Wallstrom \cite{rota}, when $\Phi_{r_1}=\cdots=\Phi_{r_n}=\Phi$,
we simply write $
\Phi^{\otimes n}_\pi$ for
$ (\Phi\otimes\cdots\otimes\Phi)_\pi$ and $\St^{[n]}_\pi$
for the corresponding measure given in (\ref{def-ito}).
Since ${C} _{\ge\widehat0}={C}$, then $\Phi^{\otimes n}_{\widehat
0}=\Phi
^{\otimes n}$, that is the product measure.
The
measure $\St^{(r_1,\ldots,r_n)}_{\widehat0}$ is called \textit{the
It{\^o} multiple stochastic measure} relative to
$\Phi_{r_1},\ldots,\Phi_{r_n}$.

As the ordinary multiple It{\^o} integral,
the It{\^o} multiple stochastic measure gives zero mass to every
diagonal set different from $C_{\widehat0}$:
\begin{proposition}
\label{0diagonals} Let $\pi\in\Pi_n$ such that $\pi>\widehat0$.
For every $C\in\mathcal{S}^{\otimes n}$, we have
\[
\St^{(r_1,\ldots,r_n)}_{\widehat0}(C_\pi)=0\qquad \mbox{a.s.}
\]
\end{proposition}
\begin{pf}
From (\ref{incompatibilitat}) we have
$(C_\pi)_{\widehat0}=\varnothing$.
\end{pf}

The basic result of Rota and Wallstrom \cite{rota}, Proposition 1,
is transferred to this situation:
\begin{proposition}
\label{strato-ito1}
%
%
\begin{equation}
\label{mesura-prod}
(\Phi_{r_1}\otimes\cdots\otimes\Phi_{r_n} )_\pi=\sum_{\sigma\ge
\pi
}\St^{(r_1,\ldots,r_n)}_\sigma
\end{equation}
and
%
%
\begin{equation}
\label{mesura-ito}
\St^{(r_1,\ldots,r_n)}_\pi=\sum_{\sigma\ge\pi}\mu(\pi,\sigma)
(\Phi
_{r_1}\otimes\cdots\otimes\Phi_{r_n} )_\sigma,
\end{equation}
where $\mu(\pi,\sigma)$ is the M\"{o}bius function defined in
Section \ref{partitions}.
\end{proposition}
\begin{pf}
The equality (\ref{mesura-prod}) is deduced from (\ref
{union}) and the definitions (\ref{def-producte}) and (\ref{def-ito}).
The equality (\ref{mesura-ito}) follows from (\ref{mesura-prod}) and
the M\"{o}ebius inversion formula (\ref{moebius}).
\end{pf}

\section{\texorpdfstring{Random measures induced by a L{\'e}vy
process}{Random measures induced by a Levy process}}
\label{random}

Let $X=\{X_t,t\in[0,T] \}$ be a L\'evy process, that is, $X$
has stationary and independent increments, is continuous in
probability, is cadlag and $X_0=0$. In all the paper we assume that $X$
has moments of all orders;
however, if the interest is restricted to multiple integral up to order
$n\ge2$, then it is enough to assume
that the process has moments up to order $2n$.

Denote the L{\'e}vy measure of $X$ by $\nu$, and by $\sigma^2$ the
variance of its Gaussian
part. The existence of moments of $X_t$ of all orders implies that
$\int_{\{\vert x\vert>1\}} \vert x\vert\nu(dx)<\infty$ and
$\int_{\mathbb{R}} \vert x\vert^n \nu(dx)<\infty, \forall n\ge
2$.
Write
%
%
\begin{eqnarray}
\label{cn}
K_1&=&E[X_1],\nonumber\\[-8pt]\\[-8pt]
K_2&=& \sigma^2+\int_{\mathbb{R}} x ^2
\nu(dx) \quad\mbox{and}\quad K_n=\int_{\mathbb{R}} x ^n
\nu(dx)<\infty,\qquad
n\ge3.\nonumber
\end{eqnarray}

From now on, take $S=[0,T]$ and $\mathcal{S}=\mathcal{B}([0,T])$. The
basic random measure $\phi$ that we consider is
the measure
induced by the process $X$ itself,
defined on the intervals by
%
%
\begin{equation}
\label{phi}
\phi(]s,t])=X_t-X_s,\qquad 0\le s\le t\le T,
\end{equation}
and extended to $\mathcal{B}([0,T])$. The measure $\phi$ is an
independently scattered random measure, that is, if $A_1,\ldots,A_n\in
\mathcal{B}([0,T])$
are pairwise disjoint, then $\phi(A_1),\ldots,\phi(A_n)$ are independent.

The random measures induced by the powers of the jumps of the process,
$\Delta X_t=X_t-X_{t-}$, are also used.
Consider the \textit{variations} of the process $X$ (see Meyer \cite{meyer},
page 319)
%
%
\begin{eqnarray}
\label{variacions}
X^{(1)}_t&=&X_t, \nonumber\\
X ^{(2)}_t&=&[X,X]_t=\sum_{0<s\le t} (\Delta X_s )^2+\sigma^2t,\\
X ^{(n)}_t&=&\sum_{0<s\le t} (\Delta X_s )^n,\qquad
n \ge3. \nonumber
\end{eqnarray}
The processes
$X^{(1)},\ldots, X^{(n)},\ldots$ are L\'evy processes such that
\[
{\mathbb{E}}\bigl[X^{(n)}_t\bigr]=K_n t\qquad \forall n\ge1.
\]
So, the centered processes,
\[
Y^{(n)}_t=X^{(n)}_t-K_n t,\qquad n\ge1,
\]
are square integrable martingales, called \textit{Teugels martingales}
(see Nualart and Schoutens \cite{nua-sch}),
with predictable quadratic covariation
\[
\bigl\langle Y^{(n)},Y^{(m)}\bigr\rangle_t= K_{n+m} t,\qquad n, m\ge1.
\]
\begin{notation}
We denote by $\phi_n$ the random measure induced by $X^{(n)}$, and for
$n=1$, $\phi_1=\phi$ (we
indistinctly use both $\phi_1$ and
$\phi$).
Every $\phi_n$ is a independently scattered random measure.
For $A, B\in\mathcal{B}([0,T])$,
\[
{\mathbb{E}}[\phi_n(A)\phi_m(B)]=K_{n+m} \int_{A\cap B}dt+
K_nK_m \int_{A} dt\int_{ B}dt.
\]
\end{notation}

We stress the following property, which is the basis of all the paper,
and is a consequence of Theorem 10.1.1 by Kwapien and Woyczynski
\cite{kwa}.
\begin{theorem}
For every $r_1,\ldots,r_n\ge1$, the random measures
$\phi_{r_1}, \ldots,\phi_{r_n}$ are jointly good random measures
on $ ([0,T]^n,\mathcal{B} ([0,T]^n) )$.
\end{theorem}

\subsection{The diagonal measures}\label{sec31}

Rota and Wallstrom \cite{rota} define the diagonal measure of order $n$
of $\phi$ as the random
measure on $[0,T]$ given by
%
%
\begin{equation}
\label{diagonal}
\bolds\Delta_n(A)=\phi^{\otimes n}( A^n_{\widehat1}),\qquad A\in\mathcal
{B}([0,T]).
\end{equation}
To identify the diagonal measures is a necessary step to study the
stochastic multiple integral.
In the case of a random measure generated by a L{\'e}vy process we show
that the diagonal measures are the measures
generated by
the variations of the process.
\begin{proposition}
\label{diagonal-simple}
For every $A\in\mathcal{B}([0,T])$ and $n\ge1$,
%
%
\begin{equation}
\label{diagonal-prop}
\bolds\Delta_n(A)=\phi_n(A),
\end{equation}
where $\phi_n$ is the random measure induced by $X^{(n)}$.
\end{proposition}
\begin{pf}
Since both $\bolds\Delta_n$ and $\phi_n$ are random measures, by
Proposition \ref{monotona} it is enough to check the equality
for $A=(0,t]$. Consider an increasing sequence of equidistributed
partitions of $[0,t]$ with the mesh going to 0;
for example, take $t_k^{(m)}=tk/2^{m}$ and let
\[
\mathcal{P}_m=\bigl\{t_k^{(m)}, k=0,\ldots, 2^m\bigr\}.
\]
To shorten the notation, write $t_k$ instead of $t_k^{(m)}$.
Consider the sets
\[
A_m= (0,t_1 ]^n\cup(t_1,t_2 ]^n\cup\cdots\cup(t_{2^m-1},t ]^n.
\]
Random measures are sequentially continuous and
$A_m \searrow(0,t]^n_{\widehat1}$, when $m\to\infty$,
so we have that
\[
\bolds\Delta_n((0,t])=\lim_m \sum_{k=0}^{2^m-1} (\phi(
(t_k,t_{k+1} ]
) )^n=
\lim_m \sum_{k=0}^{2^m-1}
(X_{t_{k+1}}-X_{t_{k}} )^n
\]
in $L^2(\Omega)$. For $n=2$,
\[
\lim_m \sum_{k=0}^{2^m-1} (X_{t_{k+1}}-X_{t_{k}} )^2=[X,X]_t=\phi_2
((0,t] )\qquad \mbox{in probability},
\]
so the proposition is true in this case.
For $n>2$,
by It{\^o}'s formula,
{\setcounter{equation}{0}
\renewcommand{\theequation}{\alph{equation}}
\begin{eqnarray}
&&\sum_{k=0}^{2^m-1} (X_{t_{k+1}}-X_{t_{k}} )^n \nonumber\\
&&\qquad=
n\sum_{k=0}^{2^m-1} \int_{t_{k}}^{t_{k+1}} (X_{s-}-X_{t_{k}} )^{n-1}
\,dX_s \nonumber\\
&&\qquad\quad{}+
\frac{1}{2} n(n-1)\sum_{k=0}^{2^m-1} \int_{t_{k}}^{t_{k+1}}
(X_{s}-X_{t_{k}} )^{n-2} \,ds\nonumber\\ 
&&\qquad\quad{}+\sum_{k=0}^{2^m-1}\sum_{t_{k}<s\le t_{k+1}} [
(X_{s}-X_{t_{k}} )^{n}- ( X_{s-}-X_{t_{k}} )^{n}\nonumber\\
&&\qquad\quad\hspace*{74.5pt}{} -n (X_{s-}-X_{t_{k}} )^{n-1} (X_{s}-X_{s-} ) ]
\nonumber\\
\label{equa}
&&\qquad=n\int_0^t \Biggl(\sum_{k=0}^{2^m-1} ( X_{s-}- X_{t_{k}} )^{n-1}\mathbf{1}
_{(t_{k},t_{k+1}]}(s) \Biggr) \,dX_s\\
\label{equb}
&&\qquad\quad{}+\pmatrix{n\cr2}\int_0^t \Biggl(\sum_{k=0}^{2^m-1} ( X_{s-}- X_{t_{k}}
)^{n-2}\mathbf{1}_{(t_{k},t_{k+1}]}(s) \Biggr) \,d[X,X]_s\\
\label{equc}
&&\qquad\quad{}+\sum_{j=3}^n\sum_{k=0}^{2^m-1}\sum_{t_{k}<s\le t_{k+1}}
\pmatrix{n\cr j} ( X_{s-}- X_{t_{k}} )^{n-j} (\Delta X_s )^j.
\end{eqnarray}
For $j=3,\ldots, n-1$, the corresponding term in (\ref{equc}) is
\begin{equation}
\label{equd}
\pmatrix{n\cr j}\int_0^t \Biggl(\sum_{k=0}^{2^m-1} ( X_{s-}- X_{t_{k}}
)^{n-j}\mathbf{1}_{(t_{k},t_{k+1}]}(s) \Biggr) \,dX^{(j)}_s.
\end{equation}}

\noindent Hence, (\ref{equa}), (\ref{equb}) and (\ref{equd}) have the same structure
\[
\int_0^tH^{(m)}_s \,dZ_s,
\]
where $H^{(m)}_s=\sum_{k=0}^{2^m-1} ( X_{s-}- X_{t_{k}} )^{r}\mathbf{1}
_{(t_{k},t_{k+1}]}(s)$ is a predictable process
and $Z$ is a semimartingale. Since $X_{s-}$ is left continuous,
\[
\lim_{m} H^{(m)}_s=0\qquad \mbox{a.s.}
\]
Moreover,
\[
\bigl\vert H^{(m)}_s \bigr\vert\le C {\sup_{0\le u \le s}} \vert X _u \vert^r,
\]
and the process
$ \{\sup_{0\le u \le s} \vert X _u \vert^r, s\in[0,t] \}$ is cadlag
and adapted, and as a consequence, it is prelocally bounded
(see pages 336 and 340 in Dellacherie and Meyer~\cite{DelMey82}).
By the dominated convergence theorem for stochastic integrals
(Dellacherie and Meyer \cite{DelMey82}, Theorem 14, page 338),
\[
\lim_{m}\int_0^tH^{(m)}_s \,dZ_s=0\qquad \mbox{in probability}.
\]
Finally, for $j=n$, the term in (\ref{equc}) is $\sum_{0<s\le t} (\Delta X_s
)^n=X^{(n)}_t$,
and the proposition is proved.
\end{pf}

Diagonal measures associated to a random measure of the form $\phi
_{r_{1}}\otimes\cdots\otimes\phi_{r_n}$ are needed.
This is an extension of the previous proposition, and it is a key
result for the sequel.
\begin{theorem}
\label{diagonal-diversos}
Let $r_1,\ldots,r_n\ge1$, $n\ge2$, and $A\in\mathcal{B}([0,T])$.
Then
\[
(\phi_{r_{1}}\otimes\cdots\otimes\phi_{r_n} )(A^n_{\widehat
1})=\bolds
\Delta
_{r_1+\cdots+r_n}(A)=\phi_{r_1+\cdots+r_n}(A).
\]
\end{theorem}
\begin{pf}
As in the proof of the last proposition and with the same notation, it
suffices to prove that
for all $t>0$
\[
\lim_m \sum_{k=0}^{2^m-1}
\bigl(X_{t_{k+1}}^{(r_1)}-X_{t_{k}}^{(r_1)} \bigr) \cdots
\bigl(X_{t_{k+1}}^{(r_n)}-X_{t_{k}}^{(r_n)} \bigr)=
\phi_{r_1+\cdots+r_n} ((0,t] )
\]
in probability. This convergence follows from Proposition \ref
{diagonal-simple} by polarization.
\end{pf}

\section{The Hu--Meyer formula: Measures}
\label{measure}
The Hu--Meyer formula gi\-ves the relationship between the product
measure $\phi^{\otimes n}$
and the It{\^o} stochastic measures
$\St_{\widehat0}^{\mathbf{r}}$.
In this section we obtain this formula for measures and in the next one
we extend it
to the corresponding integrals.

The idea of Hu--Meyer formula is the following. Given $C\in\mathcal
{B}([0,T]^n)$, we can decompose
\[
C=\bigcup_{\sigma\in\Pi_n} C_\sigma.
\]
So
\[
\phi^{\otimes n}(C)=\sum_{\sigma\in\Pi_n} \phi^{\otimes
n}(C_\sigma).
\]
Next step is to express each $\phi^{\otimes n}(C_\sigma)$ as a multiple
It{\^o} stochastic measure.
For example, take $n=3, \sigma= \{\{1\},\{2,3\} \}$ and $C=A^3$.
Then,
\[
A^3_\sigma=\{(s,t,t), s,t\in A, s\ne t\},
\]
and we will prove that
\[
\phi^{\otimes3}(A^3_\sigma)=\St_{\widehat0}^{(1,2)}(A^2).
\]
That is, both the product measure and the product set on the last two
variables collapse to produce a diagonal measure,
and since $s\ne t$, we get an It{\^o} measure.
To handle in general this property, we need
some notation.

Given a partition $\sigma\in\Pi_n$ with blocks $B_1,\ldots, B_m$, we
can order
the blocks
in agreement with the minimum element of each block.
When necessary, we assume that the blocks have been ordered with that
procedure, and we simply say that
$B_1,\ldots, B_m$ are ordered.
In that situation, we write
%
%
\setcounter{equation}{12}
\begin{equation}\label{barra}
\overline\sigma=(\#B_1,\ldots,\#B_m).
\end{equation}

We start considering a set $C=A^n$, with $A\in\mathcal{B}([0,T])$, and
later we extend the Hu--Meyer formula to an arbitrary set $C\in
\mathcal
{B}([0,T]^n)$.
\begin{theorem}
\label{hu-meyer1}
Let $A\in\mathcal{B}([0,T])$. Then
%
%
\begin{equation}
\label{hu-meyer1-form}
\phi^{\otimes n}(A^n)=\sum_{\sigma\in\Pi_n} \St_{\widehat
0}^{\overline
\sigma
}(A^{\#\sigma}).
\end{equation}
\end{theorem}

To prove this theorem we need two lemmas. The first one is an
invariance-type property of product
measures under permutations. We remember some standard notation.
\begin{notation}
We denote by
$\mathfrak{G}_n$ the set of permutations of $1,\ldots, n$.
Consider $p\in\mathfrak{G}_n$.

\begin{enumerate}
\item For a partition $\sigma\in\Pi_n$ with blocks $B_1,\ldots, B_m$,
we write
$p(\sigma)$ for the partition with blocks $W_j=p(B_j)=\{p(i), i\in
B_j\}$.
Note that in general
the blocks $W_1,\ldots,W_m$ are not ordered, even when $B_1, \ldots,
B_m$ are.
\item For a vector $\mathbf{x}=(x_1,\ldots,x_n)\in{\mathbb{R}}^n$,
we write
\[
p(\mathbf{x})= \bigl(x_{p(1)},\ldots, x_{p(n)} \bigr).
\]
Given $C\subset{\mathbb{R}}^n$, we put
\[
p(C)=\{p(\mathbf{x}) \mbox{, for } \mathbf{x}\in C\}.
\]
\end{enumerate}
\end{notation}
\begin{lemma}
\label{permutacio}
Let $p\in\mathfrak{G}_n$ and $r_1,\ldots,r_n\ge1$. Then for every
$C\in
\mathcal{B}([0,T]^n)$,
%
%
\begin{equation}
\label{permut1}
\bigl(\phi_{r_{p(1)}}\otimes\cdots\otimes\phi_{r_{p(n)}} \bigr)(p(C))=
(\phi_{r_{1}}\otimes\cdots\otimes\phi_{r_n} )(C)
\end{equation}
and
%
%
\begin{equation}
\label{permut2}
\St_{\widehat0}^{p(\mathbf{r})}(p(C))=\St_{\widehat0}^{\mathbf{r}}(C).
\end{equation}
\end{lemma}
\begin{pf}
Define the vector measure
\[
\Psi(C)= \bigl(\phi_{r_{p(1)}}\otimes\cdots\otimes\phi_{r_{p(n)}} \bigr)(p(C)).
\]
For $C=A_1\times\cdots\times A_n$, we have that
\[
p(A_1\times\cdots\times A_n)=A_{p(1)}\times\cdots\times A_{p(n)},
\]
and it is clear that
\[
\Psi(C)= (\phi_{r_{1}}\otimes\cdots\otimes\phi_{r_n} )(C).
\]
Then, equality (\ref{permut1}) follows from Proposition
\ref{monotona}.

To prove (\ref{permut2}), first note that,
by definition, the It{\^o} stochastic measure satisfies
\[
\St_{\widehat0}^{\mathbf{r}}(C)=\St_{\widehat0}^{\mathbf
{r}}(C_{\widehat0}).
\]
Moreover $ (p(C) )_{\widehat0}=p (C_{\widehat0} )$.
So it suffices to prove (\ref{permut2}) for a set
$C=C_{\widehat0}$. Denote by $\mathcal{B}^n_{\widehat0}$ the $\sigma
$-algebra
trace of $\mathcal{B}([0,T]^n)$ with $[0,T]^n_{\widehat0}$,
which
is composed by all sets~$C_{\widehat0}$, with $C\in\mathcal{B}([0,T]^n)$.
This $\sigma$-algebra is generated (on $[0,T]^n_{\widehat0}$) by the family
of rectangles $A_1\times\cdots\times A_n$, with $A_1,\ldots, A_n$
pairwise disjoint.
By Proposition~\ref{monotona}, we only need to check (\ref{permut2})
for this type of rectangle,
and the property reduces to~(\ref{permut1}).
\end{pf}

The next lemma is an important step in proving Theorem \ref{hu-meyer1}.
To have an insight into its meaning, consider the
following example: let $n=4$ and
$\sigma= \{\{1\},\{2\},\{3,4\} \}$.
With a slight abuse of notation, we can write
\[
A^4_{\ge\sigma}=\{(s,t,u,u)\dvtx s,t,u\in A\}=A^2\times A^2_{\widehat
1} .
\]
By Theorem \ref{diagonal-diversos},
\begin{eqnarray*}
(\phi_{r_1}\otimes\phi_{r_2} \otimes\phi_{r_3}\otimes\phi_{r_4} )
(A^4_{\ge\sigma} )&=&
\phi_{r_1}(A)\phi_{r_2}(A) (\phi_{r_3}\otimes\phi_{r_4}
)(A^2_{\widehat
1})\\
&=&\phi_{r_1}(A)\phi_{r_2}(A)\phi_{r_3+r_4}(A).
\end{eqnarray*}
However, if you consider $\tau= \{\{1, 3\}, \{2\}, \{4\} \}$, even though
$\tau$ and $\sigma$ have the same number of blocks with 1 element and
the same number of blocks with 2 elements
(they have the same type), the computation of
$ (\phi_{r_1}\otimes\phi_{r_2} \otimes\phi_{r_3}\otimes\phi
_{r_4} )
(A^4_{\ge\tau} )$
is not so straightforward. The lemma gives such computation. Its proof
demands some combinatorial results and it is
transferred to Appendix \ref{prova-lema}.
\begin{lemma}
\label{lema-basic} Let $r_1,\ldots,r_n\ge1$, $\sigma\in\Pi_n$ with
blocks $B_1,\ldots, B_m$ (ordered), and $A\in\mathcal{B}
([0,T])$. Then
\[
(\phi_{r_{1}}\otimes\cdots\otimes\phi_{r_n} )(A^n_{\ge\sigma})=
\prod_{j=1}^m\phi_{\sum_{i\in B_j}r_i}(A).
\]
\end{lemma}
\begin{pf*}{Proof of Theorem \ref{hu-meyer1}}
By Proposition \ref{strato-ito1},
\[
\phi^{\otimes n}(A^n)=\sum_{\sigma\in\Pi_n} \St_{\sigma}^{[n]}(A^{n}).
\]
So it suffices to prove that
\[
\St_{\sigma}^{[n]}(A^{n})=\St_{\widehat0}^{\overline\sigma}(A^{\#
\sigma}).
\]
By the second statement in Proposition \ref{strato-ito1}
we have
%
%
\begin{equation}
\label{hu-form2}
\St_{\sigma}^{[n]}(A^{n})=\sum_{\pi\in[\sigma,\widehat1]}\mu
(\sigma
,\pi)\phi
^{\otimes n}_\pi(A^n)=\sum_{\pi\in[\sigma,\widehat1]}\mu(\sigma
,\pi)
\phi^{\otimes n} (A^n_{\ge\pi} ).
\end{equation}
By Lemma \ref{lema-basic},
%
%
\begin{equation}
\label{hu-form3}
\phi^{\otimes n} (A^n_{\ge\pi} )=\prod_{V\in\pi}\phi_{\# V}(A).
\end{equation}
Let $B_1,\ldots,B_m$ be the blocks of $\sigma\in\Pi_n$ (ordered)
and write
\[
\overline\sigma=(\#B_1,\ldots,\#B_m)=(s_1,\ldots,s_m).
\]
The partition
$\pi\in[\sigma,\widehat1]$, with blocks $V_1,\ldots,V_k$, induces
a unique
partition of $\pi^*\in\Pi_m$, with blocks
$W_1,\ldots,W_k$ such that
\[
V_i=\bigcup_{j\in W_i}B_j
\]
(see Proposition \ref{bijeccio} in the \hyperref[app]{Appendix}).
Hence, for $i=1,\ldots, k$,
\[
\phi_{\# V_i}(A)=\phi_{\sum_{j\in W_i}\#B_j}(A)=\phi_{\sum_{j\in
W_i}s_j}(A).
\]
Thus,
from (\ref{hu-form3}) and Lemma \ref{lema-basic},
%
%
\begin{equation}
\label{hu-form4}
\phi^{\otimes n} (A^n_{\ge\pi} )=\prod_{W_i\in\pi^*}\phi_{\sum
_{j\in
W_i}s_j}(A)=
(\phi_{s_1}\otimes\cdots\otimes\phi_{s_m} )(A^m_{\ge\pi^*}).
\end{equation}
By (\ref{hu-form2}) and (\ref{hu-form4}) using again the bijection
between $[\sigma,\widehat1]$ and $\Pi_m$ stated
in Proposition \ref{bijeccio} in the \hyperref[app]{Appendix}, and Proposition
\ref{strato-ito1}, we obtain
\begin{eqnarray*}
\St_{\sigma}^{[n]}(A^{n})&=&\sum_{\pi\in[\sigma,\widehat1]}\mu
(\sigma
,\pi)
(\phi_{s_1}\otimes\cdots\otimes\phi_{s_m} )(A^m_{\ge\pi^*})\\
&=&\sum_{\rho\in\Pi_{m}}\mu(\widehat0,\rho)
(\phi_{s_1}\otimes\cdots\otimes\phi_{s_m} )(A^m_{\ge\rho
})=\St_{\widehat
0}^{\overline\sigma}(A^{\#\sigma}).
\end{eqnarray*}
\upqed\end{pf*}

In order to extend the Hu--Meyer formula for a general set in $\mathcal
{B}([0,T]^n)$,
we use a set function to express for an arbitrary set the contraction
from $A^n$ to $A^{\# \sigma}$. That is, given
a partition
$\sigma\in\Pi_n$, with blocks $B_1,\ldots,B_m$ ordered, we want to
contract a set $C\in\mathcal{B}([0,T]^n)$ into
a set of $\mathcal{B}([0,T]^{\# \sigma})$
according to the structure of the $\sigma$-diagonal sets.
With this purpose, define the function
%
%
\begin{eqnarray}
\label{qsigma}
q_\sigma\dvtx[0,T]^{\#\sigma} &\longrightarrow&
[0,T]^n,\nonumber\\[-8pt]\\[-8pt]
(x_1,\ldots,x_m) &\longrightarrow& (y_{1},\ldots,y_{n}), \nonumber
\end{eqnarray}
where $y_i=x_j$, if $i\in B_j$. For example, if $n=4$ and $\sigma= \{\{
1\},\{2,4\},\{3\} \}$,
\[
q_\sigma(x_1,x_2,x_3)=(x_1,x_2,x_3,x_2).
\]
Note that
\[
q_\sigma^{-1}(A^n)=A^{\#\sigma}.
\]
See Appendix \ref{q-sigma} for more details.
\begin{theorem}
\label{hu-meyer2}
Let $C\in\mathcal{B}([0,T]^n)$. Then
%
%
\begin{equation}
\label{hu-meyer2-form}
\phi^{\otimes n}(C)=\sum_{\sigma\in\Pi_n} \St_{\widehat
0}^{\overline
\sigma
}(q_\sigma^{-1}(C)).
\end{equation}
\end{theorem}
\begin{pf}
We separate the proof in two
steps. In the first one, we show that
it is enough to prove the theorem for a rectangle of the form
\[
C=A_1^{r_1}\times\cdots\times A_\ell^{r_\ell},
\]
where $A_1,\ldots,A_\ell$ are pairwise disjoint. In the second step we
check formula (\ref{hu-meyer2-form})
for those rectangles.

\textit{First step.} By Proposition \ref{monotona}, it suffices to prove
the theorem for a rectangle
$A_1\times\cdots\times A_n$. Since every rectangle can be written
as a disjoint union of rectangles such that every two components are
either equal or disjoint,
we consider one of this rectangles,
$C=A_1\times\cdots\times A_n$, where for every $i,j$, $A_i=A_j$ or
$A_i\cap A_j=\varnothing$. Now we show
that the formula (\ref{hu-meyer2-form}) applied to $C$
is invariant by permutations: specifically, we see that for any
permutation $p\in\mathfrak{G}_n$
\[
\phi^{\otimes n}(p(C))=\phi^{\otimes n}(C) \quad\mbox{and}\quad \sum_{\sigma
\in\Pi_n}
\St_{\widehat0}^{\overline\sigma}(q_\sigma^{-1}(p(C)))
=\sum_{\sigma\in\Pi_n} \St_{\widehat0}^{\overline\sigma}(q_\sigma
^{-1}(C)).
\]
The first equality is deduced from (\ref{permut1}). For the second one,
applying Proposition~\ref{perm-delta}(i), we have
\[
\St_{\widehat0}^{\overline\sigma}(q_\sigma^{-1}(p(C)))=\St_{\widehat
0}^{\overline
\sigma}\bigl(p_1^{-1}\bigl(q_{p(\sigma)}^{-1}(C)\bigr)\bigr),
\]
where $p_1 \in\mathfrak{G}_{\#\sigma}$ is the permutation that gives
the correct order of the blocks of $p(\sigma)$
(see the lines before Proposition \ref{perm-delta}).
By Lemma \ref{permutacio}
\[
\St_{\widehat0}^{\overline\sigma}\bigl(p_1^{-1}\bigl(q_{p(\sigma)}^{-1}(C)\bigr)\bigr)=
\St_{\widehat0}^{p_1(\overline\sigma)}\bigl(q_{p(\sigma)}^{-1}(C)\bigr)=
\St_{\widehat0}^{\overline{p(\sigma)}}\bigl(q_{p(\sigma)}^{-1}(C)\bigr),
\]
where the last equality is due to the fact that $p_1(\overline\sigma
)=\overline{p(\sigma)}$
by the definition of $p_1$ [see (\ref{barra2})].
Finally,
\[
\sum_{\sigma\in\Pi_n} \St_{\widehat0}^{\overline\sigma}(q_\sigma
^{-1}(C))=
\sum_{\sigma\in\Pi_n} \St_{\widehat0}^{\overline{p(\sigma
)}}\bigl(q_{p(\sigma)}^{-1}(C)\bigr),
\]
because we are adding over all the set $\Pi_n=\{p(\sigma), \sigma\in
\Pi
_n\}$.

\textit{Second step.} Consider
\[
C=A_1^{r_1}\times\cdots\times A_\ell^{r_\ell}
\]
with $A_1,\ldots, A_\ell$ pairwise disjoint and $\sum_{i=1}^\ell r_i=n$.
By Theorem \ref{hu-meyer1},
\begin{eqnarray*}
\phi^{\otimes n}(A_1^{r_1}\times\cdots\times A_\ell^{r_\ell})&=&
\prod_{i=1}^\ell\phi^{\otimes r_i}(A^{r_i}_{i})
=\prod_{i=1}^\ell\sum_{\sigma_i\in\Pi_{r_i}}\St_{\widehat
0}^{\overline
\sigma_i }(A^{\# \sigma_i}_{i})\\
&=& \sum_{\sigma_1\in\Pi_{r_1},\ldots,\sigma_\ell\in\Pi_{r_\ell}}
\St_{\widehat0}^{\overline\sigma_1,\ldots,\overline\sigma_\ell}
(A^{\# \sigma_1}_{1}\times\cdots\times A^{\# \sigma_\ell}_{\ell} ),
\end{eqnarray*}
where the last equality is due to the fact that
\[
(A^{\# \sigma_1}_{1}\times\cdots\times A^{\# \sigma_\ell}_{\ell}
)_{\widehat0}=
(A^{\# \sigma_1}_{1} )_{\widehat0}\times\cdots\times( A^{\# \sigma
_\ell
}_{\ell} )_{\widehat0},
\]
and the definition of the It\^{o} measure $\St_{\widehat0}^{\overline
\sigma
_1,\ldots,\overline\sigma_\ell}$.

Let $\tau\in\Pi_n$ be the partition with blocks
\begin{eqnarray*}
F_1&=&\{1,\ldots, r_1\},\\
F_2&=&\{r_1+1,\ldots,r_1+r_2\},\\
& \vdots\\
F_\ell&=&\{r_1+\cdots+r_{\ell-1}+1,\ldots,n\}.
\end{eqnarray*}
There is a bijection between the elements $\sigma\in\Pi_n$, with
$\sigma\le\tau$, and
$(\sigma_1,\ldots,\sigma_\ell)\in\Pi_{r_1}\times\cdots\times
\Pi
_{r_\ell}$ such that
\[
\overline\sigma=(\overline\sigma_1,\ldots,\overline\sigma_\ell)
\quad\mbox{and}\quad
q_\sigma^{-1}(A_1^{r_1}\times\cdots\times A_\ell^{r_\ell})
=A^{\# \sigma_1}_{1}\times\cdots\times A^{\# \sigma_\ell}_{\ell},
\]
where we use equality (\ref{projector2}) in the \hyperref[app]{Appendix}.
Then,
\begin{eqnarray*}
\phi^{\otimes n}(A_1^{r_1}\times\cdots\times A_\ell^{r_\ell})&=&
\sum_{\sigma\in\Pi_n, \sigma\le\tau}\St_{\widehat0}^{\overline
\sigma
}\bigl(q_\sigma^{-1}(A_1^{r_1}\times\cdots\times A_\ell^{r_\ell})\bigr)\\
&=&\sum_{\sigma\in\Pi_n}\St_{\widehat0}^{\overline\sigma}\bigl(q_\sigma
^{-1}(A_1^{r_1}\times\cdots\times A_\ell^{r_\ell})\bigr),
\end{eqnarray*}
where the last equality is due to the fact that if $\sigma\not\le
\tau$,
then $q_\sigma^{-1}(A_1^{r_1}\times\cdots\times A_\ell^{r_\ell
})=\varnothing$ [see (\ref{projector2})].
\end{pf}

\section{\texorpdfstring{Multiple It{\^o} and Stratonovich integral, and the
corresponding Hu--Meyer formula}{Multiple Ito and Stratonovich integral, and the
corresponding Hu--Meyer formula}}
\label{integral}

We extend Theorem \ref{hu-meyer2} to integrals with respect to the
random measures involved. We first
define an It{\^o}-type multiple integral and an integral with respect
to the product measure.

\subsection{\texorpdfstring{Multiple It{\^o} stochastic integral}{Multiple Ito stochastic integral}}
\label{iterada-ito}

We generalize the multiple It{\^o} integral with respect to the
Brownian motion (It{\^o} \cite{ito-multiple}; see also \cite{ito})
to a multiple integral with respect to the L{\'e}vy processes
$X^{(r_1)},\ldots, X^{(r_n)}$. As we will prove, that integral can be
interpreted
as the integral with respect to the It{\^o} stochastic measure.
The ideas used to construct this integral are mainly It{\^o}'s; however,
the fact that these processes (in general) are not centered
obstructs the classical isometry property, being substituted by an inequality.

Write $L^2_n=L^2([0,T]^n,\mathcal{B}([0,T]^n),
(dt)^{\otimes n})$. Denote by $\mathcal{E}_n^{\mathrm{Ito}}$ the set
of the
so-called It{\^o}-elementary functions, having the form
\[
f(t_1,\ldots,t_n)=\sum_{i_1,\ldots,i_n=1}^m a_{i_1,\ldots
,i_n}\mathbf{1}
_{A_{i_1}\times\cdots\times A_{i_n}}(t_1,\ldots,t_n),
\]
where $A_1,\ldots, A_m\in\mathcal{B}([0,T])$ are pairwise disjoint, and
$a_{i_1,\ldots,i_n}$ is zero if two indices
are equal. It is well known (see It{\^o} \cite{ito-multiple}) that
$\mathcal{E}_n^{\mathrm{Ito}}$ is dense in $L^2_n$.
Consider $f\in\mathcal{E}_n^{\mathrm{Ito}}$ and define the multiple
It{\^o}
integral of $f$ with respect to $X^{(r_1)},\ldots, X^{(r_n)}$ by
\[
I^{(r_1,\ldots,r_n)}_n(f)=\sum_{i_1,\ldots,i_n=1}^m a_{i_1,\ldots
,i_n}\phi
_{r_1}(A_{i_1})\cdots\phi_{r_n}(A_{i_n}).
\]
\begin{lemma} Let $f\in\mathcal{E}_n^{\mathrm{Ito}}$ and $\mathbf
{r}=(r_1,\ldots
,r_n)$. Then
\[
{\mathbb{E}}[(I_n^{\mathbf{r}}(f))^2 ]\le\alpha_{\mathbf{r}}\int
_{[0,T]^n}f^2(t_1,\ldots
,t_n) \,dt_1\cdots dt_n,
\]
where $\alpha_{\mathbf{r}}$ is a constant that depends on $r_1,\ldots,r_n$
but not on $f$.
\end{lemma}
\begin{pf}
The proof follows exactly the same steps as that of Theorem 4.1 in
Engel \cite{engel}. The key point is that
the measures $\phi_{r_i}$ can be written as
\[
\phi_{r_i}(A)=\overline{\phi}_{r_i}(A)+K_{r_i} \int_A dt,
\]
where $\overline{\phi}_{r_i}$ is the centered and independently
scattered random measure corresponding to $Y^{(r_i)}$.
\end{pf}

The extension of the multiple It{\^o} stochastic integral to $L^2_n$,
stated below, is proved as in the Brownian case (see It\^{o}
\cite{ito-multiple}).
\begin{theorem}
\label{ito-ext-th} The map
\begin{eqnarray*}
I_n^{\mathbf{r}}\dvtx\mathcal{E}_n^{\mathrm{Ito}}&\longrightarrow&
L^2(\Omega),\\
f &\longrightarrow& I_n^{\mathbf{r}}(f)
\end{eqnarray*}
can be extended to a unique linear continuous map from $L^2_n$ to
$L^2(\Omega)$. In particular,
$I_n^{\mathbf{r}}(f)$ satisfies the inequality
%
%
\begin{equation}
\label{desigualtat-integral}
{\mathbb{E}}[(I_n^{\mathbf{r}}(f))^2 ]\le\alpha_{\mathbf{r}}\int
_{[0,T]^n}f^2(t_{1},\ldots
,t_{n}) \,dt_1\cdots dt_n.
\end{equation}
\end{theorem}

As in the Brownian case, it is useful to express the multiple integral
in terms of iterated integrals of the form
\[
\int_0^T \biggl(\int_0^{t_{i_1}-}\cdots\biggl(\int
_0^{t_{i_{n-1}}-}f(t_{i_1},\ldots
,t_{i_n})\, dX^{(r_{i_n})}_{t_{i_n}} \biggr)
\cdots dX_{t_{i_2}}^{(r_{i_2})}
\biggr)\,dX^{(r_{i_1})}_{t_{i_1}},
\]
where ${i_1},\ldots,{i_n}$ is a permutation of $1,\ldots, n$.
This integral is properly defined for $f\in L^2_n$. This can be checked
using the decomposition of $X^{(r_i)}$ as a special semimartingale
$X_t^{(r_i)}=K_{r_i} t+Y^{(r_i)}_t$,
where, as we said in Section \ref{random}, $Y^{(r_i)}$ is a square
integrable martingale with
predictable quadratic variation $\langle Y^{(r_i)},Y^{(r_i)} \rangle
_t=K_{2r_i}t$. The previous
iterated integral then reduces to
a linear combination of iterated integrals of type
\[
\int_0^T \biggl(\int_0^{t_{i_1}-}\cdots\biggl(\int_0^{t_{i_{n-1}}-}f(t_1,\ldots
,t_n) \,dZ^{(n)}_{t_{i_n}} \biggr)
\cdots dZ^{(2)}_{t_{i_2}} \biggr)
\,dZ^{(1)}_{t_{i_1}},
\]
being $Z^{(j)}_t$ either $t$ or $Y_t^{(r_j)}$. Hence, at each
iteration, the integrability condition
\[
{\mathbb{E}}\biggl[\int_0^t g^2\, d\bigl\langle Z^{(i)}, Z^{(i)}\bigr\rangle\biggr]<\infty
\]
of a predictable process $g$ with respect to $Z^{(i)}$ can be easily verified.

Next proposition gives the precise expression of the multiple integral
as a sum
of iterated integrals. Since we are integrating with respect to
different processes,
we need to separate the space $[0,T]^n$ into simplexs.
\begin{proposition} Let $f\in L^2_n$. Then
\[
I^{(r_1,\ldots,r_n)}_n(f)=\sum_{p\in\mathfrak{G}_n }
\int\cdots\int_{p(\Sigma_n)} f(t_1,\ldots,t_n) \,dX^{(r_1)}_{t_1}
\cdots
dX^{(r_n)}_{t_n},
\]
where $\Sigma_n=\{0<t_1<\cdots<t_n<T\}$, and the integrals on the
right-hand side are interpreted as iterated integrals.
\end{proposition}
\begin{pf}
By linearity and density arguments, it suffices to consider a function
\[
f=\mathbf{1}_{A_1\times\cdots\times A_n},
\]
where $A_i=(s_i,t_{i}]$ are pairwise disjoint,
and a computation gives the result.
\end{pf}

When $r_1=\cdots=r_n=1$, we write $I_n(f)$ instead of $I_n^{(1,\ldots
,1)}(f)$; in that case, the multiple
It{\^o} integral enjoys nicer properties.
\begin{proposition}
\label{nicer}
\begin{enumerate}
\item Let $f\in L^2_n$. Then
\[
I_n(f)=I_n (\widetilde f ),
\]
where $\widetilde f$ is the symmetrization of $f$
%
%
\begin{equation}
\label{simetrizacio}
\widetilde f=\frac{1}{n!}\sum_{p\in\mathfrak{G}_n} f\comp p.
\end{equation}

\item Assume ${\mathbb{E}}[X_t]=0$. For $f,g \in L^2_n$,
\[
{\mathbb{E}}[I_n(f)I_m(g)]=\delta_{n,m} K_2^n n!\int
_{[0,T]^n}\widetilde{f} \widetilde g \,d\mathbf{t},
\]
where $\delta_{n,m}=1, \mbox{ if } n=m$, and 0 otherwise.

\item Let $f\in L^2_n$ be a symmetric function. Then
\[
I_n(f)=n!\int_0^T \biggl(\int_0^{t_1-}\cdots\biggl(\int
_0^{t_{n-1}-}f(t_1,\ldots
,t_n) \,dX_{t_n} \biggr) \cdots dX_{t_2} \biggr)\,dX_{t_1}.
\]
\end{enumerate}
\end{proposition}

We now state the relationship between the It{\^o} stochastic measure
$\St_{\widehat0}^{\mathbf{r}}$
and the It{\^o} multiple integral $I_n^{\mathbf{r}}$.
\begin{proposition}
\label{stochastic-int}
Let $C\in\mathcal{B}([0,T]^n)$ and $\mathbf{r}=(r_1,\ldots,r_n)$. Then
%
%
\begin{equation}
\label{ito-int-mes}
\St_{\widehat0}^{\mathbf{r}}(C)=I_n^{\mathbf{r}}(\mathbf{1}_{C}).
\end{equation}
\end{proposition}
\begin{pf}
By (\ref{desigualtat-integral}) the map $C \mapsto I_n^{\mathbf
{r}}(\mathbf{1}_{C})$ defines
a vector measure on $\mathcal{B}([0,T]^n)$.
On the left-hand side of (\ref{ito-int-mes}), the
It\^{o} measure satisfies
%
%
\begin{equation}
\label{f-0}
\St_{\widehat0}^{\mathbf{r}}(C)=\St_{\widehat0}^{\mathbf
{r}}(C_{\widehat0}).
\end{equation}
Now, look at right-hand side of (\ref{ito-int-mes}). For $\pi\in\Pi_n,$
we have that $C_\pi=C\cap[0,T]^n_\pi$. For all $\pi>\widehat0$,
\[
{\mathbb{E}}[ (I_n^{\mathbf{r}}(\mathbf{1}_{C_\pi}) )^2 ]\le
\alpha_n\int_{[0,T]^n} \mathbf{1}_{C_\pi}\,dt_1\cdots dt_n\le\alpha
_n\int
_{[0,T]^n} \mathbf{1}_{[0,T]^n_\pi} \,dt_1\cdots dt_n=0.
\]
Hence,
%
%
\begin{equation}
\label{f-00}
I_n^{\mathbf{r}}(\mathbf{1}_{C})=I_n^{\mathbf{r}}(\mathbf
{1}_{C_{\widehat0}}).
\end{equation}

From (\ref{f-0}) and (\ref{f-00}),
it suffices to prove (\ref{ito-int-mes}) for a set
$C=C_{\widehat0}$. As in the proof of the second part of Lemma
\ref{permutacio},
this can be reduced to check that equality for
a rectangle $(s_1,t_1]\times\cdots\times(s_n,t_n]$, with the intervals
pairwise disjoint.
This follows from the fact that
both sides of (\ref{ito-int-mes}) are equal to $\phi
_{r_1}((s_1,t_1])\cdots\phi_{r_n}((s_n,t_n])$.
\end{pf}

The property $I_n(f)=I_n (\widetilde f )$ is lost when the integrators
are different. However,
from Proposition \ref{stochastic-int} and (\ref{permut2}) we can deduce
the following useful property:
\begin{proposition}
\label{perm-int}
Let $f\in L^2_n$, and $\mathbf{r}=(r_1,\ldots,r_n)$, where
$r_1,\ldots,r_n\ge
1$. Consider $p\in\mathfrak{G}_n$. Then
\[
I_n^{\mathbf{r}}(f)=I_n^{p(\mathbf{r})} (f \comp p^{-1} ).
\]
\end{proposition}

\subsection{Multiple Stratonovich integral and Hu--Meyer
formula}\label{sec52}

Given a map $f\dvtx[0, T]^n\to{\mathbb{R}}$,
the integral with respect to the product measure $\phi^{\otimes n}$
is called
the multiple Stratonovich integral, and denoted by $I_n^S(f)$. Its
basic property is that the integral of
a product function
factorizes
\[
I^S_n(g_1\otimes\cdots\otimes g_n)=I^S_1(g_1)\cdots I^S_1(g_n),
\]
where
\[
I^S_1(g)=I_1(g)=\int_0^T g(t) \,dX_t.
\]
In order to construct this integral, we consider ordinary simple
functions of the measurable space
$([0,T]^n,\mathcal{B}([0,T]^n))$. Specifically,
denote by $\mathcal{E}_n^{\mathrm{Strato}}$ the set of functions with
the form
\[
f=\sum_{i=1}^k a_{i}\mathbf{1}_{C_i},
\]
where $C_i\in\mathcal{B}([0,T]^n), i=1,\ldots, k$.
For such $f$, define the multiple Stratono\-vich integral by
\[
I_n^S(f)=\sum_{i=1}^k a_{i}\phi^{\otimes n}(C_i).
\]
The integral of a simple function does not depend on its
representation, and it is linear. Moreover,
\begin{proposition}
\label{hu-meyer-funct-prop}
Let $f\in\mathcal{E}_n^{\mathrm{Strato}}$. Then we have the
Hu--Meyer formula
%
%
\begin{equation}
\label{hu-meyer-funct1}
I^S_n(f)=\sum_{\sigma\in\Pi_n} I_{\#\sigma}^{\overline\sigma
}(f\comp
q_\sigma),
\end{equation}
where the function $q_\sigma\dvtx[0,T]^{\#\sigma} \to[0,T]^n$ is
introduced in (\ref{qsigma}),
$\overline\sigma= (\#B_1,\ldots,\break\#B_m)$ is the vector whose
components are the sizes of the ordered blocks of $\sigma$,
and $I^{(s_1,\ldots, s_m)}_m$ is the multiple It{\^o} integral of order
$m$ with respect to the measures
$\phi_{s_1},\ldots, \phi_{s_m}$.
\end{proposition}
\begin{pf}
By linearity, it suffices to consider $f=\mathbf{1}_{C}$, where
$C\in\mathcal{B}([0,T]^n)$.
A generic term on the right-hand side of (\ref{hu-meyer-funct1}) is
\[
I_{\#\sigma}^{\overline\sigma}(\mathbf{1}_C\comp q_\sigma)
\]
and
\[
\mathbf{1}_C\comp q_\sigma=\mathbf{1}_{q_\sigma^{-1}(C)}.
\]
Hence, by Proposition \ref{stochastic-int},
\[
I_{\#\sigma}^{\overline\sigma}(\mathbf{1}_C\comp q_\sigma
)=\St_{\widehat
0}^{\overline
\sigma}(q_\sigma^{-1}(C)),
\]
and (\ref{hu-meyer-funct1}) follows from Theorem \ref{hu-meyer2}.
\end{pf}

Let $\sigma\in\Pi_n$, with $\#\sigma=m$, and denote by $\lambda
_\sigma
$ the image measure of the Lebesgue measure
$(dt)^{\otimes m}$ by
the function $q_\sigma\dvtx[0,T]^{m} \to[0,T]^n$. The image measure
theorem implies that for
$f\dvtx[0,T]^n\to{\mathbb{R}}$ measurable, positive
or $\lambda_\sigma$-integrable,
%
%
\begin{equation}
\label{image}
\int_{[0,T]^n} f(t_1,\ldots,t_n) \,d\lambda_\sigma(t_1,\ldots,t_n)=
\int_{[0,T]^m} f(q_\sigma(t_1,\ldots,t_m)) \,dt_1 \cdots dt_m.\hspace*{-28pt}
\end{equation}
Define on $\mathcal{B}([0,T]^n)$ the measure
\[
\Lambda_n=\sum_{\sigma\in\Pi_n}\lambda_\sigma,
\]
and write $L^2(\Lambda_n)$ for $L^2([0,T]^n,\mathcal
{B}([0,T]^n),\Lambda_n)$.

In order to extend the multiple Stratonovich integral we need the
following inequality of norms:
\begin{lemma}
Let $f\in\mathcal{E}_n^{\mathrm{Strato}}$. Then
%
%
\begin{equation}
\label{desigualtat-strato}
{\mathbb{E}}[ (I_n^S(f) )^2 ]\le C\int_{[0,T]^n}f^2\, d\Lambda_n,
\end{equation}
where $C$ is a constant.
\end{lemma}
\begin{pf}
By
(\ref{hu-meyer-funct1}), (\ref{desigualtat-integral}) and (\ref{image}),
\begin{eqnarray*}
{\mathbb{E}}[ (I_n^S(f) )^2 ]& \le& C \sum_{\sigma\in\Pi_n}
{\mathbb{E}}\bigl[ \bigl(I_{\#\sigma}^{\overline\sigma}(f\comp q_\sigma) \bigr)^2
\bigr]\\
&\le& C \sum_{\sigma\in\Pi_n}\int_{[0,T]^{\#\sigma}}(f\comp
q_\sigma)^2\,
dt_1\cdots dt_{\#\sigma}\\
&=& C \sum_{\sigma\in\Pi_n}\int_{[0,T]^{n}}f^2 \,d\lambda_\sigma\\
&=&C\int
_{[0,T]^{n}}f^2 \,d\Lambda_n.
\end{eqnarray*}
\upqed\end{pf}

The main result of the paper is the following theorem:
\begin{theorem}
The map $I_n^S\dvtx\mathcal{E}_n^{\mathrm{Strato}}\to L^2(\Omega)$
can be extended to a unique linear continuous map from $L^2(\Lambda_n)$
to $L^2(\Omega)$,
and we have
the Hu--Meyer formula
%
%
\begin{equation}
\label{hu-meyer-funct}
I^S_n(f)=\sum_{\sigma\in\Pi_n} I_{\#\sigma}^{\overline\sigma
}(f\comp
q_\sigma).
\end{equation}
\end{theorem}
\begin{pf}
The extension of $I_n^S$ to a continuous map on $L^2(\Lambda_n)$ is
proved using a density argument and inequality (\ref{desigualtat-strato}).
To prove the Hu--Meyer formula,
let $f\in L^2(\Lambda_n)$ and $\{f_k, k\ge1\}\subset\mathcal
{E}_n^{\mathrm{Strato}}$ such that $\lim_kf_k=f$ in
$L^2(\Lambda_n)$. For every $\sigma\in\Pi_n$, we have
$\lim_k f_k\comp q_\sigma=f\comp q_\sigma$ in $L^2_{\# \sigma}$; hence,
from Theorem \ref{ito-ext-th} the It{\^o} integrals on the right-hand
side of
(\ref{hu-meyer-funct})
converge, and the formula follows from Proposition
\ref{hu-meyer-funct-prop}.
\end{pf}
\begin{remarks}
\begin{enumerate}[(1)]
\item[(1)] Let $g_1,\ldots,g_n\in{L}^{2n}([0,T],dt)$. Then
$g_1\otimes\cdots
\otimes g_n\in L^2_n(\Lambda_n)$ and
\[
I^S_n(g_1\otimes\cdots\otimes g_n)=I^S_1(g_1)\cdots I^S_1(g_n).
\]
This result is easily checked for simple functions $g_1,\ldots,g_n$ and
extended to the
general case by a density argument.

\item[(2)] In order to prove the Hu--Meyer formula for $I_n^S$ it is
enough to
assume that the process $X$ has moments up to order $2n$.

\item[(3)] For $\sigma\in\Pi_n, \sigma>\widehat0$, the measure
$\lambda_\sigma$
is singular with respect
to the\vspace*{1pt} Lebesgue measure
on $[0,T]^n$. For example, for $n=2$ and $\sigma=\widehat1$,
let $D=\{(t,t), t\in[0,T]\}$ be the diagonal of $[0,T]^2$. Then
$\lambda_{\widehat1}$ is concentrated in $D$, that has
zero Lebesgue measure, but $\lambda_{\widehat1}$ is nonzero
\[
\lambda_{\widehat1}(D)=\int_{[0,T]^2}\mathbf{1}_{D}(s,t) \,d\lambda
_{\widehat
1}(s,t)=\int
_{[0,T]}\mathbf{1}_{D}(t,t) \,dt=T.
\]

\item[(4)] As in the Brownian case (see
\cite{johnson,sole-tra,hu-meyer2,masani} and the references therein),
there are other procedures
to construct the multiple Stratonovich
integral. The main difficulty in every approach
is that the usual condition $f\in L^2_n$ in It{\^o}'s theory is not
sufficient to guarantee
the multiple Stratonovich integrability of $f$.
The reason is that one needs
to control the behavior of $f$ on the diagonal sets $[0,T]_\sigma^n$
that have zero Lebesgue
measure when $\sigma>\widehat0$.
We solve this difficulty using the norm induced by
the measure $\Lambda_n$, which seems to be appropriate for dealing
with the diagonal sets,
avoiding in this way the difficulty of a manageable definition of the
\textit{traces}.
\end{enumerate}
\end{remarks}

When
the function $f\in L^2(\Lambda_n)$ is symmetric, the Hu--Meyer formula
can be considerably simplified.
We show that we can assume that symmetry on $f$ without loss of generality.
\begin{proposition} Let $f\in L^2(\Lambda_n)$. Then
$I^S_n(f)=I^S_n(\widetilde
f)$, where $\widetilde f$ is the symmetrization of $f$
[see (\ref{simetrizacio})].
\end{proposition}
\begin{pf}
The proof is straightforward for $f=\mathbf{1}_C$, $C\in\mathcal
{B}([0,T]^n)$,
using Lemma \ref{permutacio}.
By linearity the equality $I_n^S(f)=I_n^S(\widetilde f)$
is extended to $\mathcal{E}^{\mathrm{Strato}}_n$,
and by density to $L^2(\Lambda_n)$.
\end{pf}

Next we show the Hu--Meyer formula for a symmetric function $f$.
In general (for $f$ symmetric), the function $f\comp q_\sigma$ is
nonsymmetric, but as we will see
in the proof of the next theorem,
its multiple It{\^o} integral depends only
on the block structure of $\sigma$ (the type of $\sigma$).
For example, with $n=3$, $f(t_1,t_2,t_3)=t_1t_2t_3$ and $\sigma= \{\{
1\}
, \{2,3\} \}$, we have that
\[
f(q_\sigma(t_1,t_2))=t_1 t_2^2,
\]
that is nonsymmetric. Its integral is
\[
I^{\overline{\sigma}}_{\#\sigma}(f\comp q_{\sigma
})=I^{(1,2)}_2(f\comp
q_\sigma)= I^{(1,2)}_2(t_1 t_2^2).
\]
Take $\pi= \{\{1,3\}, \{2\} \}$. Then $f(q_\pi(t_1,t_2))=t_1^2 t_2$ and
\[
I^{\overline{\pi}}_{\#\pi}(f\comp q_{\pi})=
I^{(2,1)}_2(t_1^2 t_2)=I^{(1,2)}_2(t_1 t_2^2),
\]
where the last equality is due to Proposition
\ref{perm-int}.

We use the following notation: given nonnegative integers $r_1,\ldots
,r_k$ such that $\sum_{i=1}^kir_i=n$,
we write
\[
[r_1, r_2,\ldots, r_k]=(\underbrace{1,\ldots,1}_{r_1},\underbrace
{2,\ldots
,2}_{r_2},\ldots).
\]
Note that this corresponds to $\overline\sigma$ when
\[
\sigma= \bigl\{\{1\}, \ldots,\{r_1\}, \{r_1+1,r_1+2\},\ldots,\{
r_1+2r_2-1,r_1+2r_2\},\ldots\bigr\}.
\]
We also write $q_{r_1,\ldots,r_k}$ for $q_\sigma$, with $\sigma$ the
above partition.
\begin{theorem} Let $f\in L^2_n(\Lambda_n)$ be a symmetric function. Then
%
%
\begin{equation}
\label{hu-meyer-simetric}
I^S_n(f)=\sum\frac{n!}{r_1!(2!)^{r_2}r_2!\cdots
(k!)^{r_k}r_k!}I_{r_1+\cdots+r_k}^{[r_1,\ldots,r_k]}
(f\comp q_{r_1,\ldots,r_k} ),
\end{equation}
where the sum is extended over all nonnegative integers $r_1,\ldots,r_k$
such that $\sum_{i=1}^kir_i=n$,
for $k=1,\ldots,n$.
\end{theorem}
\begin{pf}
Let $f\in L^2(\Lambda_n)$ symmetric.
For every $\sigma\in\Pi_n$ and $p\in\mathfrak{G}_n$,
\begin{eqnarray*}
\label{inv-perm}
I_{\#p(\sigma)}^{\overline{ p(\sigma)}}\bigl(f\comp q_{p(\sigma)}\bigr)&=_{\fontsize{8.36}{10.36}\selectfont{\mbox{(\ref{equa})}}}&
I_{\#p(\sigma)}^{\overline{ p(\sigma)}}(f\comp p^{-1}\comp q_{\sigma
}\comp p_1^{-1})\\
&=_{\fontsize{8.36}{10.36}\selectfont{\mbox{(\ref{equb})}}}&I_{\#p(\sigma)}^{\overline{ p(\sigma)}}(f\comp q_{\sigma
}\comp p_1^{-1})
=_{\fontsize{8.36}{10.36}\selectfont{\mbox{(\ref{equc})}}}I_{\#\sigma}^{\overline{\sigma}}(f\comp q_{\sigma}),
\end{eqnarray*}
where (\ref{equa}) is due to Proposition \ref{perm-delta}(ii), the equality
(\ref{equb}) follows from the symmetry of $f$ and
(\ref{equc}) from Proposition \ref{perm-int} and the fact that $p_1$ gives the
correct order of $p(\sigma)$ [see (\ref{barra2})].
This implies that all the partitions that have the same number of
blocks of 1 element, the
same number with two elements, etc. (i.e., they have the same type)
give the same It{\^o} multiple integral
in the Hu--Meyer formula.
To obtain (\ref{hu-meyer-simetric}) it suffices to count the number of
partitions of $\{1,\ldots,n\}$ with
$r_1$ blocks with 1 element, $r_2$ blocks with 2 elements$,\ldots,r_k$
blocks with $k$ elements, which is
\[
\frac{n!}{r_1!(2!)^{r_2}r_2!\cdots(k!)^{r_k}r_k!}.
\]
\upqed\end{pf}

\textit{Final remark.}
One may expect that by decomposing the L\'evy process into a sum of two
independent processes, one with the small jumps and the
other with the large ones, the assumption of the existence of moments
could be avoided.
However, this decomposition
introduces dramatic changes to the context of the work, and such an
extension is beyond the scope and purposes of the present paper.

\section{Special cases}\label{sec6}
\subsection{Brownian motion}\label{sec61}
When $X=W$ is a standard Brownian motion,
\[
\phi_2([0,t])=t \quad\mbox{and}\quad \phi_n=0,\qquad n\ge3.\vadjust{\goodbreak}
\]
It follows that in the Hu--Meyer formula only the partitions with all
blocks of cardinality 1 or 2 give a contribution,
and all the It{\^o} integrals are a mixture of multiple stochastic
Brownian integrals and Lebesgue integrals.
We can organize the sum according the number of blocks of two elements.
For a partition having
$j$ blocks of $2$ elements, and $f\in L^2_n(\Lambda_n)$ symmetric, the
multiple It{\^o} integral is
\begin{eqnarray*}
&&I^{[n-2j,j]}_{n-j}(f)\\
&&\qquad=\int_{[0,T]^{n-j}}f(s_1,\ldots,s_{n-2j},\\
&&\hspace*{84.4pt}t_1,t_1,\ldots
,t_{j},t_{j}) \,dW_{s_1}\cdots dW_{s_{n-2j}} \,d t_{1}\cdots
d t_{j}\\
&&\qquad=I_{n-2j} \biggl(\int_{[0,T]^{j}}f(\bolds\cdot, t_1,t_1,\ldots
,t_{j},t_{j}) \,d t_{1}\cdots
d t_{j} \biggr),
\end{eqnarray*}
where the last equality is due to a Fubini-type theorem. Therefore,
%
%
\begin{eqnarray}
\label{hu-meyer-brow}
I_n^S(f)
&=&\sum_{j=0}^{[n/2]}\frac{n!}{(n-2j)!j!2^j}\nonumber\\[-8pt]\\[-8pt]
&&\hspace*{17.6pt}{}\times I_{n-2j}
\biggl(\int_{[0,T]^{j}}f(\bolds\cdot,t_1,t_1,t_2,t_2,\ldots,t_{j},t_{j})\,
d t_{1}\cdots
d t_{j} \biggr),\nonumber
\end{eqnarray}
which is the classical Hu--Meyer formula (see \cite{hu-meyer1}).

On the other hand, in the measure $\Lambda_n$ only participate the measures
$\lambda_\sigma$ corresponding to the partitions above mentioned.
Consider the measure $\ell_2=\lambda_{\widehat1}$ on $[0,T]^2$,
that is, for a positive or $\ell_2$ integrable function $h$,
\[
\int_{[0,T]^2}h(s,t) \,d\ell_2(s,t)=\int_{[0,T]}h(t,t)\,dt.
\]
Given the partition $\sigma\in\Pi_n$,
\[
\sigma= \bigl\{\{1\},\ldots,\{n-2j\},\{n-2j+1,n-2j+2\},\ldots, \{n-1,n\}
\bigr\},
\]
we have
\[
\lambda_\sigma=(dt)^{\otimes(n-2j)}\otimes\ell_2^{\otimes j}.
\]

\subsection{Poisson process}\label{sec62}
Let $N_t$ be a standard Poisson process with intensity 1, and consider
the process $X_t=N_t-t$.
For every $n\ge2$,
\[
X^{(n)}_t=N_t=X_t+t,
\]
and hence, a multiple It{\^o} integral can be reduced to a linear
combination of multiple integrals where all
the integrators are $dX$ or $dt$. For $f\in L^2_n(\Lambda_n)$
symmetric, each integral
$I_{r_1+\cdots+r_k}^{[r_1,\ldots,r_k]}
(f\comp q_{r_1,\ldots,r_k})$ in (\ref{hu-meyer-simetric}) can be
expressed in terms of
the number of Lebesgue integrals that appear
\begin{eqnarray*}
\hspace*{-5pt}&&I_{r_1+\cdots+r_k}^{[r_1,\ldots,r_k]}
(f\comp q_{r_1,\ldots,r_k})\\
\hspace*{-5pt}&&\qquad=\sum_{j=0}^{r_2+\cdots+r_k}\hspace*{-0.9pt}I_{r_1+\cdots+r_k-j} \Biggl(
\int_{[0,T]^j} \Biggl(\mathop{\sum_{l_1,\ldots,l_j=r_1+1}}_{\mathrm
{different}\mbox{ }
}^{r_1+\cdots+r_k}\hspace*{-0.9pt}
(f\comp q_{r_1,\ldots,r_k})\\
\hspace*{-5pt}&&\hspace*{204.5pt}{}\times (t_1,\ldots,t_{r_1+\cdots+r_k}) \Biggr)\,dt_{l_1}\cdots dt_{l_j} \Biggr),
\end{eqnarray*}
and the Hu--Meyer formula of Sol{\'e} and Utzet \cite{sole-pois1} can
be deduced from this expression.

\subsection{Gamma process and subordinators}\label{sec63}

A subordinator is a L\'{e}vy process with increasing paths. An
important example of a subordinator with moments of all orders is the
Gamma process,
denoted by $\{G_t, t\ge0\}$, which is the L\'{e}vy process
corresponding to an exponential law of parameter 1. Its
L\'{e}vy
measure is
\[
\nu(dx)=\frac{e^{-x}}{x} \mathbf{1}_{\{x>0\}}(x) \,dx.
\]
The law of $G_t$ is Gamma with mean $t$ and scale parameter equal to
one. A Gamma process can be represented as the sum of its jumps, that are
all positive,
\[
G_t=\sum_{0<s\le t}\Delta G_s.
\]
The L\'{e}vy measure of $G^{(n)}$ is (see Schoutens \cite{Sch00})
\[
\nu_n(dx)=\frac{e^{-x^{1/n}}}{nx}\mathbf{1}_{\{x>0\}}(x) \,dx,\qquad n\ge1,
\]
and the Teugels martingales are
\[
Y^{(n)}_t=\sum_{\{0< s \le t\}} (\Delta G_s )^n-(n-1)! t,\qquad n\ge1.
\]

In this case, unlike the Brownian motion and the Poisson process, the
Hu--Meyer formula does not simplify, due to the fact that the diagonal
measures cannot
be expressed in a simple way in terms of, say, the process and
a deterministic measure. However, for a Gamma process, and in general,
for a subordinator without drift (see below for the definition) with
moments of all orders,
both the multiple It\^{o} and Stratonovich integrals can be computed
pathwise integrating with respect to an ordinary measure.
This is a multivariate
extension of the property that states that the stochastic integral and
the pathwise Lebesgue--Stieljes integral with respect to a semimartingale
of bounded variation are equal; such property was proved for the
integral with respect to a L\'{e}vy process of bounded variation by
Millar \cite{Mil70}
under weak conditions, and part of our proof follows his scheme.

Let $X=\{X_t, t\ge0\}$ be a subordinator. The L\'{e}vy--It\^{o}
representation of $X$ takes the form
\[
X_t=\gamma_0t+ \sum_{0<s\le t}\Delta X_s
\]
with $\gamma_0\ge0$ (see Sato \cite{Sat99}, Theorems 21.5 and 19.3).
The number $\gamma_0$ is called the \textit{drift} of the subordinator,
and we
will assume that $\gamma_0=0$.
Consider the sequence of stopping times $\{T_k, k\ge1\}$ with disjoint
graphs that exhaust
the jumps
of $X\dvtx\Delta X_{T_k}\neq0, \forall k\ge1$, and $X$ only has jumps
on these times (see, e.g., Dellacherie and Meyer \cite{DelMey82},
Theorem B, page XIII,
for a construction of
this sequence). Denote by $J_n$ the set of $n$-tuples $(T_{i_1},\ldots,
T_{i_n})$, with $T_{i_j}\le T$, and all entries different.
For $r_1,\ldots,r_n\ge1$, define a measure on $[0,T]^n$ by
\[
m_{r_1,\ldots,r_n}=\sum_{(T_{i_1},\ldots,T_{i_n})\in J_n} (\Delta
X_{T_{i_1}} )^{r_1}\cdots(\Delta X_{T_{i_n}} )^{r_n}
\delta_{(T_{i_1},\ldots, T_{i_n})},
\]
where $\delta_a$ is a Dirac measure at point $a$, with the convention
that the sum is 0 if $J_n=\varnothing$.
We have the following property:
\begin{proposition}
Let $X=\{X_t, t\ge0\}$ be a subordinator without drift and with
moments of all orders. With the preceding notation, for every $f\in L^2_n$,
%
%
\begin{equation}
\label{path}
I_n^{(r_1,\ldots,r_n)}(f)=\int_{[0,T]^n}f \,d m_{r_1,\ldots,r_n}\qquad
\mbox{a.s.}
\end{equation}
\end{proposition}
\begin{pf}
First, note two facts:

(a) $m_{r_1,\ldots,r_n}$ is a finite measure
\begin{eqnarray*}
m_{r_1,\ldots,r_n}([0,T]^n)&=&\sum_{(T_{i_1},\ldots, T_{i_n})\in J_n}
(\Delta X_{T_{i_1}} )^{r_1}\cdots(\Delta X_{T_{i_n}} )^{r_n}\\
&\le&\sum_{T_{i_1}\le T,\ldots, T_{i_n}\le T} (\Delta X_{T_{i_1}}
)^{r_1}\cdots(\Delta X_{T_{i_n}} )^{r_n}\\
&=&X_T^{(r_1)}\cdots X_T^{(r_n)}<\infty.
\end{eqnarray*}

(b) If the intervals $(s_1,t_1],\ldots,(s_n,t_n]$ are pairwise disjoint,
then (\ref{path}) is true for $f=\mathbf{1}_{(s_1,t_1]\times\cdots
\times
(s_n,t_n]}$. The proof is straightforward.

We separate the proof of the proposition in two steps.

\textit{Step} 1. Formula (\ref{path}) is true for every map $f\dvtx
[0,T]^n\to
{\mathbb{R}}$ $\mathcal{B}_0$-measurable and bounded, where $\mathcal{B}_0$
is the $\sigma$-field on $[0,T]^n$ generated by the rectangles
$(s_1,t_1]\times\cdots\times(s_n,t_n]$,
with
$(s_1,t_1],\ldots,(s_n,t_n]$ pairwise disjoint.

To prove this claim we use a convenient monotone class theorem.
Denote by $\mathcal{H}$ the family of functions that satisfy (\ref
{path}); it is a vector space such that:
\begin{longlist}
\item$1\in\mathcal{H}$.
\item If $f_m\in\mathcal{H}$, $0\le f_m\le K$ for some constant $K$,
and $f_m \nearrow f$, then $f\in\mathcal{H}$.
\end{longlist}
To see (i), consider the dyadic partition of $[0,T]$ with mesh
$2^{-k}$, write
\[
B_j= \bigl((j-1)T2^{-k},jT2^{-k} \bigr],\qquad j=1,\ldots, 2^k,
\]
and define
\[
f_k=\mathop{\sum_{j_1,\ldots,j_n}}_{\mathrm{different}}
\mathbf{1}_{B_{j_1}\times\cdots\times B_{j_n}}.
\]
By the remark (b) at the beginning of the proof,
\[
I_n^{(r_1,\ldots,r_n)}(f_k)=\int_{[0,T]^n}f_k \,d m_{r_1,\ldots,r_n}.
\]
Moreover, $f_k \nearrow1$ out off the diagonal sets $[0,T]^n_\sigma$,
with $\sigma\ne\widehat0$, and
then $f_k \nearrow1$ a.e. with respect to the Lebesgue measure, and in $L^2_n$.
Therefore,
\[
\lim_{k}I_n^{(r_1,\ldots,r_n)}(f_k)=I_n^{(r_1,\ldots,r_n)}(1).
\]
On the other hand, for every $\omega$, the measure $m_{r_1,\ldots,r_n}$
does not charge on any of the above mentioned diagonal sets. Thus, the
convergence
$f_k \nearrow1$ is also $m_{r_1,\ldots,r_n}$-a.e. By the monotone
convergence theorem,
\[
\lim_k \int_{[0,T]^n}f_k \,d m_{r_1,\ldots,r_n}=\int_{[0,T]^n}1\, d
m_{r_1,\ldots,r_n},
\]
and (i) follows.

Point (ii) is deduced directly from the monotone convergence theorem
and taking into account that under the conditions in (ii)
we have $f_m \to f$ in $L^2_n$.

Again by remark (b) above, the indicator of a set $(s_1,t_1]\times
\cdots\times(s_n,t_n]$, with
$(s_1,t_1],\ldots,(s_n,t_n]$ pairwise disjoint, is in $\mathcal{H}$, and
this family of sets is closed by intersection.
By the monotone class theorem, it follows that all bounded $\mathcal
{B}_0$-measurable functions are in $\mathcal{H}$.

\textit{Step} 2. Extension of (\ref{path}) to all $f\in L^2_n$. First,
note that
$\mathcal{B}_0$ is the $\sigma$-field generated by the Borelian sets
$B\in\mathcal{B}([0,T]^n)$ such that $B\subset[0,T]^n_{\widehat0}$.
Then,\vspace*{-1pt} given $B\in\mathcal{B}([0,T]^n)$, the indicator $\mathbf
{1}_{B\cap
[0,T]^n_{\widehat0}}$
is $\mathcal{B}_0$ measurable. Let $f\in L_n^2$, and assume $f\ge0$.
There is a sequence of simple (and then bounded) functions
such that $0\le f_m \nearrow f$.
Define\vspace*{-1pt} $f^{0}_m=f_m \mathbf{1}_{[0,T]^n_{\widehat0}}$, which is
$\mathcal{B}_0$
measurable, and
$f_m^0 \nearrow f$ a.e. with respect to the Lebesgue measure. The
convergence is also in $L^2_n$, and then
$\lim_m I_n(f_m^0)=I(f)$. On the other hand, $f_m^0 \nearrow f$,
$m_{r_1,\ldots,r_n}$-a.e. so
\[
\lim_m\int_{[0,T]^n}f_m^0 \,d m_{r_1,\ldots,r_n}= \int_{[0,T]^n}f
\,d m_{r_1,\ldots,r_n}.
\]
By Step 1, we get the result. For a general
$f\in L^2_n$, decompose $f=f^+-f^-$.
\end{pf}

Finally, for a subordinator without drift and with moments of all
orders, the multiple Stratonovich measure can
be identified with the $n$-fold product measure of $\phi=\sum_{k}
\Delta X_{T_k} \delta_{T_k}$. So
for $f\in\mathcal{E}_n^{\mathrm{Srato}}$, by definition,
\[
I_n^S(f)=\int_{[0,T]^n}f \,d\phi^{\otimes n}.
\]
Using similar arguments as in the previous proposition, but easier, it
is proved that
\[
I_n^S(f)=\int_{[0,T]^n}f \,d\phi^{\otimes n}\qquad \forall f\in L^2(\Lambda_n).
\]
Then, the Hu--Meyer formula can be transferred to a pathwise context.

\begin{appendix}\label{app}
\section*{Appendix}

\subsection{The isomorphism $[\sigma,\widehat1]\simeq\Pi_{\#\sigma}$}

Fix a partition $\sigma\in\Pi_n$, with blocks $B_1,\ldots, B_m$.
Let $\pi
\ge\sigma$, with blocks
$V_1,\ldots, V_k$; each block $V_i$ is the union of some of the blocks
$B_1,\ldots, B_m$. Hence, we
can consider the partition $\pi^*\in\Pi_m$ that gives the relationship
between the $V_i$'s and the $B_j$'s, that is,
$\pi^*$ has blocks $W_1,\ldots,W_k$ defined by
\[
V_i=\bigcup_{j\in W_i }B_j,\qquad i=1,\ldots,k.
\]
\setcounter{proposition}{0}
\begin{proposition}
\label{bijeccio}
Let $\sigma\in\Pi_n$ with $\#\sigma=m$. With the above notation,
the map
\begin{eqnarray*}
[\sigma,  \widehat1]&\longrightarrow& \Pi_m,\\
\pi&\mapsto&\pi^*
\end{eqnarray*}
is a bijection and, for $\pi,\tau\in[\sigma,\widehat1]$,
\[
\pi\le\tau\quad\Longleftrightarrow\quad\pi^*\le\tau^*.
\]
Moreover,
\[
\mu^{(n)}(\sigma,\pi)=\mu^{(m)}(\widehat0,\pi^*),
\]
where $\mu^{(r)}$ is the M{\^o}bius function on $\Pi_r$.
\end{proposition}

The proof is straightforward.

\subsection{Permutations and partitions}
\label{a-permutations}

Let $p\dvtx\{1,\ldots,n\}\longrightarrow\{1,\ldots, n\}$ be a permutation.
This application induces
a bijection on $\Pi_n$, and a bijection on ${\mathbb{R}}^n$.
Specifically:

1. For a subset $B\subset\{1,\ldots,n\}$ we denote by $p(B)$ the
image of $B$ by $p$
\[
p(B)=\{p(j)\mbox{, for } j\in B\}.
\]
Given a partition $\sigma\in\Pi_n$, with blocks $B_1,\ldots,B_m$,
let $p(\sigma)$ be the partition with blocks $W_1,\ldots, W_m$ defined
by $W_j=p(B_j)$. Note that in general
the blocks $W_1,\ldots,W_m$ are not ordered. The application
\begin{eqnarray*}
p\dvtx\Pi_n&\longrightarrow&\Pi_n,\\
\sigma&\mapsto& p(\sigma)
\end{eqnarray*}
is a bijection and for $\sigma,\tau\in\Pi_n$,
\[
\sigma\le\tau\quad\Longleftrightarrow\quad p(\sigma)\le p(\tau).
\]
This last property is clear, because if $V\in\tau$, and
$V=B_{r_1}\cup
\cdots\cup
B_{r_k}$,
then
\[
p(V)=p(B_{r_1})\cup\cdots\cup p(B_{r_k}).
\]

Further, this application is compatible with the relationship
introduced in Section \ref{subsets}
%
%
\begin{equation}
\label{equivalencia}
i\sim_\sigma j \quad\Longleftrightarrow \quad p(i) \sim_{p(\sigma)} p(j).
\end{equation}

2. For a vector $\mathbf{x}=(x_1,\ldots,x_n)\in{\mathbb{R}}^n$, we write
\[
p(\mathbf{x})= \bigl(x_{p(1)},\ldots, x_{p(n)} \bigr),
\]
and the application $\mathbf{x}\mapsto p(\mathbf{x})$ determines a
bijection on
${\mathbb{R}}^n$, that we also denote by~$p$.
For a set $C\subset{\mathbb{R}}^n$, we write
\[
p(C)=\{p(\mathbf{x})\mbox{, for } \mathbf{x}\in C\}.
\]
In particular,
for $A_1,\ldots, A_n\subset{\mathbb{R}}$,
\[
p(A_1\times\cdots\times A_n)=A_{p(1)}\times\cdots\times A_{p(n)}.
\]
Notice that if we look for the position of a particular set, say $A_1$,
in $p(A_1\times\cdots\times A_n)$, we find it at place $p^{-1}(1)$
\begin{eqnarray*}
&A_{i_1}\times\cdots\times\fbox{A$_1$}\times\cdots\times
A_{i_n}&\\
&\uparrow&\\
&p^{-1}(1).&
\end{eqnarray*}
This last observation gives some light to the next property:
\begin{proposition}
\label{permutacio-gran}
Consider $p\in\mathfrak{G}_n$, $C\subset{\mathbb{R}}^n$ and $\sigma
\in\Pi_n$.
\begin{longlist}
\item$p(C_\sigma)= (p(C) )_{p^{-1}(\sigma)}$.

\item$p(C_{\ge\sigma})= (p(C) )_{\ge p^{-1}(\sigma)}$. In particular,
$p(A^n_{\ge\sigma})=A^n_{\ge p^{-1}(\sigma)}$.
\end{longlist}
\end{proposition}
\begin{pf}
(i) Let $\mathbf{x}=(x_1,\ldots, x_n)\in p(C_\sigma)
\subset p(C)$.
Write
\[
p^{-1}(\mathbf{x})=\bigl(x_{p^{-1}(1)},\ldots,x_{p^{-1}(n)}\bigr)=(y_1,\ldots
,y_n)=\mathbf{y}
\in C_\sigma.
\]
Therefore,
\[
\mathbf{y}\in C\quad \mbox{and}\quad y_i=y_j \quad\Longleftrightarrow\quad i\sim_\sigma j.
\]
The condition on the right is
equivalent to $p^{-1}(i) \sim_{p^{-1}(\sigma)} p^{-1}(j)$ [see (\ref
{equivalencia})].
So, returning to the $\mathbf{x}$'s,
\[
\mathbf{x}\in p(C)\quad \mbox{and}\quad x_{p^{-1}(i)}=x_{p^{-1}(j)}
\quad\Longleftrightarrow\quad p^{-1}(i) \sim_{p^{-1}(\sigma)} p^{-1}(j).
\]
Call $p^{-1}(i)=r$ and $p^{-1}(j)=s$. We have
\[
\mathbf{x}\in p(C)\quad \mbox{and}\quad x_{r}=x_{s}\quad \Longleftrightarrow\quad r \sim
_{p^{-1}(\sigma)} s.
\]
Hence, $\mathbf{x}\in(p(C) )_{p^{-1}(\sigma)}$.

The reciprocal inclusion is analogous.

(ii) Applying (i),
\begin{eqnarray*}
p (C_{\ge\sigma} )&=&p \biggl(\bigcup_{\pi\ge\sigma} C_\pi\biggr)=\bigcup
_{\pi\ge
\sigma}
p (C_\pi)
=\bigcup_{\pi\ge\sigma} (p(C) )_{p^{-1}(\pi)}\\
&=&\bigcup_{\tau\ge p^{-1}(\sigma)} (p(C) )_{\tau}=
(p(C) )_{\ge p^{-1}(\sigma)}.
\end{eqnarray*}
\upqed\end{pf}

Consider a partition $\sigma\in\Pi_n$ with blocks $B_1,\ldots, B_m$
(ordered). If the elements of each block
are consecutive numbers,
then,
%
%
\begin{equation}
\label{descomposicio}
A^n_{\ge\sigma}=\bigtimes{} A^{\#B_j}_{\widehat1}.
\end{equation}
When $\sigma$ does not fulfill the previous condition, the expression
(\ref{descomposicio}) is not valid.
However, since we are interested in computing $ (\phi_{r_{1}}\otimes
\cdots\otimes\phi_{r_n} )(A^n_{\ge\sigma})$,
thanks to
Lemma \ref{permutacio}, we fortunately can permute both the set and the
product measure
to make things work.
The next proposition is essential for this purpose.
\begin{proposition}
\label{particio} Let $A\subset{\mathbb{R}}$ and $\sigma\in\Pi_n$
be a
partition, with blocks $B_1,\ldots,B_m$ (ordered).
There is a permutation
$p\in\mathfrak{G}_n$ such that
\[
p (A^n_{\ge\sigma} )=\bigtimes A^{\#B_j}_{\widehat1}.
\]
\end{proposition}
\begin{pf}
Write $s_j=\#B_j, j=1,\ldots,n$, and let $p'\in\mathfrak
{G}_n$ such that
\begin{eqnarray*}
p'(B_1)&=&\{1,\ldots,s_1\},\\
p'(B_2)&=&\{s_1+1,\ldots,s_1+s_2\}\\
\vdots&
\end{eqnarray*}
Take $p=(p')^{-1}$ and apply Proposition \ref{permutacio-gran}(ii).
\end{pf}

\subsection{\texorpdfstring{Proof of Lemma \protect\ref{lema-basic}}{Proof of Lemma
4.4}}
\label{prova-lema}
We prove the following lemma:\vspace*{10pt}

\textsc{Lemma \ref{lema-basic}.}\quad
Let $r_1,\ldots,r_n\ge1$, $\sigma\in\Pi_n$ with blocks $B_1,\ldots,
B_m$ (ordered), and $A\in\mathcal{B}[0,T]$. Then
\[
(\phi_{r_{1}}\otimes\cdots\otimes\phi_{r_n} )(A^n_{\ge\sigma})=
\prod_{j=1}^m\phi_{\sum_{i\in B_j}r_i}(A).
\]
\begin{pf}
Let $B_1,\ldots,B_m$ be the blocks of $\sigma$ ordered. If
$\sigma$ is such that
$A^n_{\ge\sigma}=\bigtimes A^{\# B_j}_{\widehat1}$,
by Theorem \ref{diagonal-diversos},
\[
(\phi_{r_{1}}\otimes\cdots\otimes\phi_{r_n} )(A^n_{\ge\sigma})
=\prod_{j=1}^{m} \biggl(\bigotimes_{i\in B_j}\phi_{r_i} (A^{\# B_j}_{\widehat
1} ) \biggr)=
\prod_{j=1}^{m}\phi_{\sum_{i\in B_j}r_i}(A).
\]
For the general case,
let $p$ be the permutation given by Proposition \ref{particio} and
write $V_j=p^{-1}(B_j), j=1,\ldots, m$.
By Proposition \ref{particio} (first), and $\#B_j=\#V_j$ (second), we have
\[
p (A^n_{\ge\sigma} )=\bigtimes A^{\#B_j}_{\widehat1}=\bigtimes
A^{\#V_j}_{\widehat1}.
\]
By Lemma \ref{permutacio} and the first part of the proof,
\begin{eqnarray*}
(\phi_{r_{1}}\otimes\cdots\otimes\phi_{r_n} )(A^n_{\ge\sigma})
&=& (\phi_{r_{p(1)}}\otimes\cdots\otimes\phi_{r_{p(n)}} )(p(A^n_{\ge
\sigma}))\\
&=& (\phi_{u_1}\otimes\cdots\otimes\phi_{u_n} )\Biggl(\bigtimes A^{\#
V_j}_{\widehat1}\Biggr)\\
&=&\prod_{j=1}^m\phi_{\sum_{i\in V_j}u_i}(A),
\end{eqnarray*}
where $u_i=r_{p(i)}$.
For every $j=1,\ldots, m$,
\[
\sum_{i\in V_j}u_i=\sum_{i\in V_j}r_{p(i)}=\sum_{i\in
p^{-1}(B_j)}r_{p(i)}=\sum_{i\in B_j}r_{i}.
\]
\upqed\end{pf}

\subsection[The expansion function]{The function $q_\sigma$}
\label{q-sigma}

Given a partition
$\sigma\in\Pi_n$, with blocks $B_1,\ldots, B_m$ (ordered), the
function $q_\sigma$ [see (\ref{qsigma})] is defined by
\begin{eqnarray*}
q_\sigma\dvtx[0,T]^m &\longrightarrow& [0,T]^n,\\
(x_1,\ldots,x_m) &\to& (y_{1},\ldots,y_{n}),
\end{eqnarray*}
where $y_i=x_j$, if $i\in B_j$. This function is a bijection between
$[0,T]^m$ and $[0,T]^n_{\ge\sigma}$, and it
is Borel measurable because
%
%
\begin{equation}
\label{projector2}
q^{-1}_\sigma(A_1\times\cdots\times A_n)= \biggl(\bigcap_{i\in B_1} A_i
\biggr)\times\cdots\times\biggl(\bigcap_{i\in B_m}
A_i \biggr).
\end{equation}

Given a partition $\sigma\in\Pi_n$, with blocks (ordered)
$B_1,\ldots
,B_m$ and a permutation $p\in\mathfrak{G}_n$, as we commented,
the blocks of $p(\sigma)$ in general are not ordered. It is convenient
to consider the
permutation $p_1 \in\mathfrak{G}_m$ that gives the correct order of
the blocks of $p(\sigma)$, that means,
$p_1(1)=i$ if $p(B_i)$ is the first block of $p(\sigma)$, $p_1(2)=j$
if $p(B_j)$ is the second block,
and so on; in other words,
\[
p\bigl(B_{p_1(1)}\bigr),\ldots, p\bigl(B_{p_1(m)}\bigr)
\]
are the blocks of $p(\sigma)$ ordered.
Remember that we defined [see (\ref{barra})] the $m$-dimensional vector
$\overline\sigma=(\#B_1,\ldots,\#B_m)$. Then
%
%
\begin{equation}\label{barra2}
p_1(\overline\sigma)=\overline{p (\sigma)}.
\end{equation}
\begin{proposition}
\label{perm-delta}
Consider $\sigma\in\Pi_n$, with $\#\sigma=m$, $p\in\mathfrak
{G}_n$, and
let $p_1 \in\mathfrak{G}_m$
be the permutation that gives the correct order of the blocks of
$p(\sigma)$.
\begin{longlist}
\item For $A_1,\ldots, A_n\in\mathcal{B}([0,T])$,
\[
q_\sigma^{-1} \bigl(p(A_1\times\cdots\times A_n) \bigr)=p_1^{-1} \bigl(q_{p(\sigma
)}^{-1}(A_1\times\cdots\times A_n) \bigr).
\]
\item$ p^{-1}\comp q_\sigma=q_{p(\sigma)}\comp p_1$.
\end{longlist}
\end{proposition}
\begin{pf}
(i) Let $B_1,\ldots, B_m$ the blocks of $\sigma$ (ordered). We have
\begin{eqnarray*}
q_\sigma^{-1} \bigl(p(A_1\times\cdots\times A_n) \bigr)&=&
q_\sigma^{-1} \bigl(A_{p(1)}\times\cdots\times A_{p(n)}\bigr)\\
&=& \biggl(\bigcap_{i\in B_1} A_{p(i)} \biggr)\times\cdots\times\biggl(\bigcap_{i\in
B_m} A_{p(i)} \biggr)\\
&=& \biggl(\bigcap_{i\in p(B_1)} A_{i} \biggr)\times\cdots\times\biggl(\bigcap_{i\in
p(B_m)} A_{i} \biggr)\\
&=&G_1\times\cdots\times G_m,
\end{eqnarray*}
where
\[
G_j=\bigcap_{i\in p(B_j)} A_{i},\qquad j=1,\ldots,m.
\]
Since $p_1$ gives the correct order of $p(B_1),\ldots, p(B_m)$,
\begin{eqnarray*}
q_{p(\sigma)}^{-1}(A_1\times\cdots\times A_n)&=&
\biggl(\bigcap_{i\in p(B_{p_1(1)})} A_{i} \biggr)\times\cdots\times\biggl(\bigcap
_{i\in p(B_{p_1(m)})} A_{i} \biggr)\\
&=&G_{p_1(1)}\times\cdots\times G_{p_1(m)}=p_1(G_1\times\cdots\times G_m),
\end{eqnarray*}
and then
\[
p_1^{-1} \bigl(q_{p(\sigma)}^{-1}(A_1\times\cdots\times A_n) \bigr)=G_1\times
\cdots\times G_m.
\]

(ii) Consider $\mathbf{y}=(y_1,\ldots,y_n) \in[0,T]^n$.
Since $\{\mathbf{y}\}=\{y_1\}\times\cdots\times\{y_n\}$,
\[
(p^{-1}\comp q_\sigma)^{-1}(\{y\})=q_\sigma^{-1}(\{p(\mathbf{y})\}
)=_{(*)}p_1^{-1}\bigl(q^{-1}_{p(\sigma)}(\{\mathbf{y}\})\bigr)
=\bigl(q_{p(\sigma)}\comp p_1\bigr)^{-1}(\{y\}),
\]
where the equality (*) is due to part (i).
\end{pf}
\end{appendix}

%

%
\printaddresses

\end{document}